\documentclass[11pt]{article} 

\usepackage[utf8]{inputenc} 

\usepackage{cancel}
\usepackage{float}

\usepackage{geometry} 
\usepackage{graphicx} 
\usepackage[english]{babel}

\usepackage{booktabs} 
\usepackage{array} 
\usepackage{paralist} 
\usepackage{verbatim} 
\usepackage{subfig} 
\usepackage{amsfonts}
\usepackage{amssymb}
\usepackage{amsmath}
\usepackage{amsthm}
\usepackage{makecell}

\usepackage{epigraph}
\usepackage{bbm}

\usepackage{hyperref}
\hypersetup{
    colorlinks=true,
    citecolor=blue
}

\renewcommand{\leq}{\leqslant}
\renewcommand{\geq}{\geqslant}

\newcommand {\matl}{\left[ \begin{matrix}}
\newcommand {\matr}{\end{matrix}\right]}

\newcommand {\EE}{ \mathbb E }

\newcommand {\Var}{\mathbf {Var}}

\newcommand{\cP}{\mathcal{P}}

\DeclareMathAlphabet{\mathbbmsl}{U}{bbm}{m}{sl}

\usepackage[round]{natbib}
\usepackage{wrapfig}
\usepackage{subfig}
\usepackage[capitalize,noabbrev]{cleveref}
\usepackage{authblk}
\usepackage{changepage}

\usepackage{setspace}

\usepackage{macros}

\title{Randomized and exchangeable improvements of\\ 
Markov's, Chebyshev's and Chernoff's inequalities}
\author[1,2]{Aaditya Ramdas}
\author[2]{Tudor Manole}
\affil[1]{Machine Learning Department, Carnegie Mellon University} 
\affil[2]{Department of Statistics and Data Science, Carnegie Mellon University}
\affil[ ]{\texttt{ \{aramdas,tmanole\}@andrew.cmu.edu  }}

\date{\today}

\newtheorem{theorem}{Theorem}[section]

\newtheorem{proposition}[theorem]{Proposition}

\newtheorem{corollary}{Corollary}[theorem]

\newtheorem{remark}[theorem]{Remark}

\begin{document}

\maketitle

\begin{abstract}
We present simple randomized and exchangeable improvements of Markov's inequality, as well as Chebyshev's inequality and Chernoff bounds.  Our variants are never worse and typically strictly more powerful than the original inequalities. The proofs are short and elementary, and can easily yield similarly randomized or exchangeable versions of a host of other inequalities that employ Markov's inequality as an intermediate step.
We point out some simple statistical applications involving tests that combine dependent e-values. In particular, we uniformly improve the power of universal inference, and obtain tighter betting-based nonparametric confidence intervals. Simulations reveal nontrivial gains in power (and no losses) in a variety of settings.
\end{abstract}

\section{Introduction}

Consider a standard probability space $(\Omega, \mathcal{F}, P)$.
Markov's inequality (MI) states that for any real-valued random variable $X$ defined on this space and constant $a>0$, we have
\begin{equation}\label{eq:M}
P(|X| \geq 1/a) \leq a \cdot \EE[|X|].
\end{equation}
If $X$ is nonintegrable, then the inequality trivially holds, so the reader may implicitly assume going forward that we deal with integrable $X$, without loss of generality. We now present three inequalities that are all strictly stronger than Markov's inequality, with an eye towards applications and improvements of other inequalities.

\subsection{The exchangeable Markov inequality (EMI)}

In contrast, the 
following stronger
version of Markov's inequality
was   recently noted by~\cite{manole2021sequential}. Given its seemingly basic and fundamental nature, it may exist elsewhere in the literature. 
We will refer to it as 
the {\it exchangeable Markov inequality (EMI)}.
\begin{theorem}[Exchangeable Markov Inequality]
\label{thm:emi}
Let $X_1,X_2,\dots$ form an exchangeable sequence of integrable random variables, meaning that the joint distribution of $(X_1,\dots,X_n)$ equals that of $(X_{\pi(1)},\dots,X_{\pi(n)})$ for any $n\geq 1$ and any permutation $\pi$ of $\{1,\dots, n\}$. Then, for any $a > 0$,
\begin{equation}\label{eq:ex-M}
P\left( \exists t \geq 1: \left| \frac1t \sum_{i=1}^t X_i \right| \geq 1/a\right) \leq P\left( \exists t \geq 1:  \frac1t \sum_{i=1}^t \left|X_i \right| \geq 1/a\right)  \leq a \cdot \EE[|X|].
\end{equation}
\end{theorem} 
This is clearly a strictly stronger statement than Markov's inequality which effectively makes an identical claim only at $t=1$ (or when all random variables are equal almost surely). 

The proof of Theorem~\ref{thm:emi} is short (the first inequality is obvious, so we focus on the second). The exchangeability of $X_1,X_2,\dots$ implies that the process $\left(\sum_{i=1}^n \left|X_i \right|/n\right)_{n \geq 1}$ forms a nonnegative reverse martingale; indeed, this follows
from the fact that it can be rewritten as $\big(\bbE\big[|X_1|\big|\calE_n\big]\big)_{n\geq 1}$,
where $(\calE_n)_{n\geq 1}$ denotes the exchangeable 
filtration generated by $X_1, X_2, \dots$  (defined for instance in~\cite{manole2021sequential} and references therein). The claim then follows by invoking the time-reversed Ville inequality (recapped in Theorem~\ref{thm:reverse_ville} below, for completeness).

We point out immediately that the above inequality results in an improvement to subsampled universal inference~\citep{wasserman2020universal}, where repeated sample splitting gives rise to exchangeable ``split likelihood ratios'' on different subsets of the data. We return to this topic in more detail in Section~\ref{sec:univ}.

A useful corollary of the EMI is as follows. Let $X_1,\dots,X_n$ be any set of (potentially nonexchangeable) arbitrarily dependent random variables. Let $\pi$ be a uniformly random permutation of $\{1,\dots,n\}$. Then, for any $a>0$, 
\begin{equation}\label{eq:ex-M2}
P\left( \sup_{1 \leq t \leq n}  \frac1t \sum_{i=1}^t \left|X_{\pi(i)} \right| \geq 1/a\right)  \leq a \cdot \frac{\EE[|X_1|+|X_2|+\dots+|X_n|]}{n}.
\end{equation}
Here, the original random variables are effectively made exchangeable by the random permutation, thus allowing us to invoke the original inequality. 

We state a final variant of the inequality. Suppose we take $N$ arbitrarily dependent random variables and put them in a bag. Suppose $X_{\pi(1)},\dots,X_{\pi(n)}$ be $n$ samples drawn uniformly at random with or without replacement from this bag.  Then, we have
\begin{equation}\label{eq:ex-M3}
P\left( \sup_{1 \leq t \leq n}  \frac1t \sum_{i=1}^t \left|X_{\pi(i)} \right| \geq 1/a\right)  \leq a \cdot \frac{\EE[|X_1|+|X_2|+\dots+|X_N|]}{N}.
\end{equation}
This holds because the sampling process induces the exchangeability required for~\eqref{eq:ex-M} to be invoked on the otherwise non-exchangeable random variables.

The aforementioned three variants of EMI are all relatively weak, in the sense that they do not really improve with increasing $t,n,N$: there is no concentration of measure really happening, and indeed there cannot really be any since we have assumed so little about the underlying random variables (a first moment for each, and either an exchangeable or an arbitrary dependence structure). Hence we are under no illusions that these can be ``much'' stronger than the original Markov's inequality. 

\subsection{Uniformly-randomized Markov inequality (UMI)}

Another improvement of Markov's inequality is the following \emph{uniformly-randomized Markov inequality (UMI)}.

\begin{theorem}[Uniformly-randomized Markov Inequality]
\label{thm:umi}
Let $X$ be a nonnegative random variable, 
and $U\sim \text{Unif}(0,1)$ an independent  random variable. Then, for any $a > 0$,
\begin{equation}\label{eq:UMI}
P(X \geq U/a) = \EE[\min(a X,1)] \leq a \cdot \EE[X].
\end{equation}
The equality above is nontrivial even for nonintegrable $X$. 
Further, that equality becomes an inequality $\leq$ if $U$ is stochastically larger than uniform. Last, if $X$ is bounded, that is $X \in [0,C]$ almost surely for some $C > 0$, then the inequality in~\eqref{eq:UMI} holds with equality for any $a \leq 1/C$.
\end{theorem}
The proof is simple:
\begin{equation}\label{eq:UMI-proof}
P(X \geq U/a) = \EE[P( U \leq a X \mid X)] = \EE[\min(a X,1)],
\end{equation}
yielding the claim.
As mentioned above, if $X \leq 1/a$ almost surely, then~\eqref{eq:UMI} holds with equality, removing all looseness in Markov's inequality. Further, the statement is nontrivial for nonintegrable $X$. For example, when $X$ is a standard Cauchy distribution and $a=0.05$, we obtain $P(|X| \geq U/a)$ actually \emph{equals} $(20\pi + \log(401) - 40 \tan^{-1}(20))/(20\pi) \approx 0.127$, but Markov's inequality only states that $P(|X| \geq 1/a)$ is \emph{at most} $0.127$.


Clearly, UMI is strictly stronger than MI. Noting that $\EE[U]=1/2$ gives an intuitive idea of the extent of the gain. The proof is so simple that the above inequality may have been previously noted by other authors. But it appears to not be commonly taught or broadly known, and we have not found any source containing it so far.
Since the UMI also holds if $U$ is stochastically larger than uniform, 
Theorem~\ref{thm:umi} also
gives an alternative proof to Markov's inequality (by choosing $U=1$).

\begin{remark} 
What happens when Markov's inequality holds with equality? This only happens in a few rare cases, but it does happen, and in these cases UMI also holds with equality. For example, take the discrete distribution with point masses at 0 and 2, with probability half each, so that the mean is one. For $a=1/2$, we have $\{X \geq 2\}=\{X=2\}$ and thus $P(X \geq 2) = 1/2$, meaning Markov's holds with equality. In this case, UMI also holds with equality, because the event $\{X \geq 2U\}$ is almost surely equal to the event $\{X=2\}$, since $U>0$ almost surely.
To summarize, Markov's inequality is only tight for a discrete random variable taking values in $\{0,1/a\}$, while the UMI holds with equality for any random variable taking values in $[0,1/a]$.
\end{remark}

We again note that the UMI
 improves the recent method of universal inference~\citep{wasserman2020universal},
by rejecting the null hypothesis when the split likelihood ratio (or its variants like the crossfit or subsampled likelihood ratio) exceed $U/\alpha$ rather than $1/\alpha$. We return to this ``randomized universal inference'' in Section~\ref{sec:univ}.

Another particular application of UMI was recently and independently noted as a passing remark, for a different context, in a recent preprint by~\cite{ign2022values}. It is known that, given an e-value $E$ (defined in Section~\ref{sec:e-val}), the mapping from $E \mapsto 1/E$ is an admissible e-to-p calibrator~\citep{vovk2021values} (that converts an e-value into a p-value). But this admissibility result only applies to deterministic maps. After proving that $P/E$ is a valid p-value where $P$ is a  p-value and $E \perp P$, \cite{ign2022values} note that this implies $U/E$ is also a valid (randomized) e-to-p calibrator that converts an e-value $E$ into a p-value  using independent randomization $U$, thus rendering $1/E$ as inadmissible if randomization is allowed. 
We expand on this type of application in Section~\ref{sec:e-val}.

Last, it turns out that the UMI and EMI can be combined into a single statement that implies both; this is postponed to Theorem~\ref{thm:umi+emi}.

\subsection{Additively-randomized Markov inequality (AMI)}\label{sec:ami}

Instead of multiplicative randomization (as done in the previous subsection), we may consider additive randomization. In this context,~\cite{huber2019halving} recently proved an interesting ``smoothed Markov'' inequality: for any nonnegative  $X$ and constant $\epsilon>0$, 
\begin{equation}\label{eq:huber}
P(X + B \geq \epsilon) \leq \frac{\EE[X]}{2\epsilon},
\end{equation}
where $B$ is an independent uniform random variable on $[-\epsilon,\epsilon]$.
However, while the right-hand side halves the bound of Markov's inequality (as his paper title suggests), the left-hand side is not directly comparable. Indeed, when calculating $P(X \geq \epsilon-B)$, the right-hand side $\epsilon - B$ is not always bigger than $\epsilon$. Hence, it appears that the bound is not in general comparable to Markov's inequality.
Huber derives many interesting consequences of this inequality, including extensions to Chebyshev's and Chernoff's inequalities, and once again they appear incomparable to the original inequalities. 

We present the following result, which (unlike Huber's) is  stronger than Markov's inequality.
\begin{proposition}[Additively-randomized Markov Inequality]
\label{prop:ami}
Let $\epsilon>0$. 
Given a nonnegative random variable $X$, 
and an independent random variable $A\sim \text{Unif}(0,\epsilon)$, it holds that
\begin{equation}\label{eq:AMI}
P(X \geq \epsilon - A) \leq \frac{\EE[X]}{\epsilon}.
\end{equation}
\end{proposition}
The proof is simple. Since $P(A \geq \epsilon - x)= x/\epsilon$ for $x \leq \epsilon$, the left-hand side simplifies to 
\[
P(A \geq \epsilon - X) = \EE[P(A \geq \epsilon - X | X)] = \EE[\min(X,\epsilon)]/\epsilon,
\]
implying the claim.

It is also not hard to see that our 
additively-randomized Markov inequality (Proposition~\ref{prop:ami}) is actually equivalent to our earlier uniformly-randomized Markov inequality
(Theorem~\ref{thm:umi}). To see this, write $A = \epsilon U$, where $U$ is uniform on $[0,1]$. Then $P(X \geq \epsilon - A) = P(X \geq \epsilon(1-U)) = P(X \geq \epsilon U')$, where $U' = 1-U$ is also uniform on $[0,1]$. Writing $a=1/\epsilon$ equates the two claims.

In Appendix~\ref{sec:huber}, we further discuss the relationship of the above bound to Huber's smoothed Markov inequality. We show, in particular,
how the two types of bounds can be used to derive each other, despite having different interpretations and implications.

In the rest of this paper, we will continue to use our multiplicative version, because we think the resulting expressions are cleaner, but readers may find the additive version more useful in some settings, which is our reasoning for recording Proposition~\ref{prop:ami} as a separate result.

\subsection{Contributions and paper outline}

Having already introduced multiple new generalizations of Markov's inequality, we next point out several new concentration inequalities that result out of their use and/or combination. Section~\ref{sec:cheby} derives randomized and exchangeable improvements of Chebyshev's inequality (with Appendix~\ref{sec:cantelli} containing an improvement to Cantelli's inequality, as a way of exemplifying a more general proof technique). Section~\ref{sec:Chernoff} does the same for Hoeffding's inequality, and points out that the same techniques improve any Chernoff bound (extensions of the Bernstein and empirical Bernstein inequalities are in Appendix~\ref{sec:bernstein}).
Section~\ref{sec:Ville} shows how to derive a randomized
improvement of Ville's inequality for forward supermartingales.
Then, Section~\ref{sec:rev-ville} randomizes the reverse Ville's inequality (for reverse submartingales), and as a consequence derives an inequality that combines the strengths of the EMI and UMI into a single inequality (Theorem~\ref{thm:rand_reverse_ville}).
After presenting some of these improvements (and omitting others for brevity), we describe some statistical applications.
Section~\ref{sec:e-val} produces more powerful tests using arbitrarily dependent e-values.  Section~\ref{sec:univ} uniformly improves different versions of the recent universal inference~\citep{wasserman2020universal} methodology. Section~\ref{sec:betting} improves betting-based tests and confidence intervals that are an exciting development in nonparametric statistics~\citep{waudby2020estimating}. Section~\ref{sec:simulations} explores the (often large) improvements in power obtained in a variety of simulations. Section~\ref{sec:sum} contains an extended discussion, including concerns about reproducibility of randomized tests, and the possibility of avoiding external randomization entirely by utilizing the internal randomness of the data. Section~\ref{sec:conclusion} contains a brief conclusion.

\section{Randomized and Exchangeable Chebyshev inequality}\label{sec:cheby}

In the rest of this section and paper, $U$ is always a uniform random variable on $[0,1]$, independent of all other random variables.

\subsection{Uniformly-randomized Chebyshev inequality}

For any random variable $X$ having variance (at most) $\sigma^2$, Chebyshev's inequality states that
\begin{equation}\label{eq:cheby-one}
P\left( | X - \EE X| \geq k \sigma \right) \leq 1/k^2.
\end{equation}
The proof is transparent: one just squares both terms within the probability on the left-hand side, and applies Markov's inequality. 
More generally, consider $n \geq 1$ i.i.d.\  random variables $X_1,\dots,X_n$ with variance $\sigma^2$, and define $\bar X_n := (X_1+\dots+X_n)/n$. Chebyshev's inequality implies that for any $k>0$,
\begin{equation}\label{eq:cheby}
P\left( | \bar X_n - \EE X| \geq k \frac{\sigma}{\sqrt{n}}\right) \leq 1/k^2.
\end{equation}
Our \emph{uniformly-randomized Chebyshev inequality} reads as follows:
\begin{theorem}[Uniformly-randomized Chebyshev Inequality]
\label{thm:rand-chebyshev}
Let $X_1, \dots, X_n$ be i.i.d.\ random variables
with variance $\sigma^2$, and let $U \sim \text{Unif}(0,1)$
be an independent random variable. Then, 
for any $k > 0$,
\begin{equation}\label{eq:cheby-u}
P\left( |\bar X_n - \EE X| \geq k \sigma  \sqrt{\frac{U}{n}}\right) \leq 1/k^2.
\end{equation}
As before, the same result holds when $U$ is stochastically larger than uniform.
\end{theorem}
This is clearly a tighter claim than the original: the probability of exceeding a random (but always smaller) threshold is identical. Noting that $\EE[\sqrt U] = 2/3$ while $\Var(\sqrt U)=1/18$ gives an intuitive idea of the extent of the gain: the obtained confidence intervals will be $2/3$ as wide using our improved inequality. 

The proof is similar to the UMI~\eqref{eq:UMI-proof}:
\begin{align*}
P\left( |\bar X_n - \EE X| \geq k \sigma  \sqrt{\frac{U}{n}}\right) &= \EE\left[P \left(U \leq \frac{n |\bar X_n - \EE X|^2}{k^2\sigma^2}   \Big| X\right) \right] 
\leq  \EE\left[ \frac{n |\bar X_n - \EE X|^2}{k^2\sigma^2} \right] = 1/k^2.
\end{align*}
The sole inequality exists, despite $U$ being exactly uniform, because $\frac{n |\bar X_n - \EE X|^2}{k^2\sigma^2}$ could be larger than one. If one replaces the resulting inequality by $\EE\left[ \frac{n |\bar X_n - \EE X|^2}{k^2\sigma^2} \wedge 1\right]$, then the inequality would turn to equality, but the subsequent equality ($=1/k^2$) would turn into an inequality. 

(It is clear that the i.i.d.\  assumption above, and in the following sections, is made for convenience, and can be weakened to give many variants of the above.)

It is also possible to prove a 
randomized version of Cantelli's inequality~\citep{cantelli1929sui}, 
which is a one-sided analogue of Chebyshev's inequality. We present such a result in Appendix~\ref{app:bentkus}.

\subsection{Exchangeable Chebyshev inequality}

If the data are exchangeable rather than i.i.d., we have the following claim:
\begin{theorem}[Exchangeable Chebyshev Inequality]
\label{thm:exch-chebyshev}
    If $X_1,X_2,\dots$ is a sequence of exchangeable random variables with variance at most $\sigma^2$, then for any $k>1$, 
    \begin{equation}\label{eq:ECI}
    P\left( \sup_{m \geq 1} | \bar X_m - \EE X| \geq k \sigma \right) \leq 1/k^2.    
    \end{equation}
    If the random variables are further assumed to be i.i.d., then for any $n\geq 1$, we have
    \begin{equation}
\label{eq:independent-chebyshev}
P\left(\sup_{m\geq n}  | \bar X_m - \EE X| \geq \frac{k\sigma}{ \sqrt n}\right)
\leq 1/k^2.
\end{equation}
\end{theorem}

Note that \eqref{eq:ECI} improves~\eqref{eq:cheby-one} by recovering it either at $m=1$ or when all random variables are identical.
Unlike the various bounds presented for the i.i.d.\ case,~\eqref{eq:ECI}
does not improve 
as the sample size increases. But improvements are not possible without further assumptions,
because when all $X_i$ are identically equal to $X$, the statement reduces to a claim about $X$, with no role for concentration of measure.

The proof is simple: defining $R_m := | \bar X_m - \EE X|^2/\sigma^2$, a short calculation invoking Jensen's inequality reveals that $(R_m)_{m\geq 1}$ is a nonnegative reverse submartingale with $\EE[R_1] \leq 1$. The time-reversed Ville inequality (Theorem~\ref{thm:reverse_ville}) then yields our claim. Essentially the same proof can be obtained by defining $Y_i = |X_i - \EE X|^2/\sigma^2$. Applying the exchangeable Markov inequality, we infer that $P(\sup_{m \geq 1} \bar Y_m \geq k^2) \leq 1/k^2$, which implies~\eqref{eq:ECI} by Jensen's inequality. For the second part, 
one may apply the same proof as above, but now
only considering the submartingale
$(R_m)_{m\geq n}$ starting at time $n$.
Since $\Var[R_n] = \sigma^2/n$ under the i.i.d.
assumption,
equation~\eqref{eq:independent-chebyshev}
follows. 

\begin{remark}
    The above inequalities are not to be confused with Kolmogorov's generalization of Chebyshev's inequality which states that  if $X_1,X_2,\dots,X_n$ are i.i.d.\ (not just exchangeable) random variables with variance at most $\sigma^2$, then for any $k>1$ and $n \geq 1$, 
\begin{equation}\label{eq:Kolmogorov}
    P\left( \sup_{m \leq n} | \sum_{i=1}^m (X_i - \EE X)| \geq k \sigma \sqrt{n} \right) \leq 1/k^2. 
    \end{equation}
    Technically, the i.i.d.\ assumption above can be weakened to just independence, and even to a martingale dependence assumption, but we omit these for simplicity.
\end{remark}

\section{Randomized and Exchangeable Chernoff bounds}\label{sec:Chernoff}

We present randomized and exchangeable variants of Chernoff bounds below.
We remark that every Chernoff bound, including matrix concentration bounds, or self-normalized concentration inequalities, are improved via the same technique presented below~\citep{howard2020time}, but it is impractical to develop every one of these, so we just pursue three of them: the Hoeffding bound below,  Bernstein's inequality in Appendix~\ref{sec:bernstein},
and the empirical Bernstein inequality
in Appendix~\ref{sec:empirical_bernstein}.

\subsection{Uniformly-randomized Hoeffding inequality}\label{sec:hoeff}

Recall that $X$ is called $\sigma$-subGaussian if for any constant $\lambda$,
\begin{equation}\label{eq:subG}
\EE[\exp(\lambda (X-\EE[X]))] \leq \exp(\lambda^2\sigma^2/2).
\end{equation}
Hoeffding's inequality states that if $X$ is $\sigma$-subGaussian, then for any $\epsilon > 0$,
\begin{equation}\label{eq:hoeff-one}
P\left( X - \EE[X] \geq \sigma \epsilon \right) \leq \exp(-\epsilon^2/2).
\end{equation}

More generally, considering $n$ i.i.d.\ $\sigma$-subGaussian random variables $X_1,\dots,X_n$, the above inequality implies that
\begin{equation}\label{eq:hoeff}
P\left( \bar X_n - \EE[X] \geq \sigma \epsilon \right) \leq \exp(-n\epsilon^2/2).
\end{equation}
Setting the right-hand side to equal $\alpha\in (0,1)$, it can be rewritten as
\begin{equation}\label{eq:hoeff2}
P\left( \bar X_n - \EE[X] \geq \sigma \sqrt{\frac{2\log(1/\alpha)}{n}} \right) \leq \alpha.
\end{equation}
An identical bound can also be derived on the other tail. As we shall soon see, the proof proceeds by  multiplying both sides within the probability by $\lambda > 0$, exponentiating both sides, then applying Markov's inequality, and finally tuning $\lambda$.

Our \emph{uniformly-randomized Hoeffding inequality} reads as follows.
\begin{theorem}[Uniformly-randomized Hoeffding Inequality]
\label{thm:rand-chernoff}
Let $X_1, \dots, X_n$ be i.i.d.\ $\sigma$-subGaussian random variables, and let $U \sim \text{Unif}(0,1)$ be an independent
random variable. Then, for all $\alpha \in (0,1)$, 
\begin{equation}\label{eq:hoeff-U}
   P\left( \bar X_n - \EE[X] \geq  \sigma \sqrt{\frac{2\log(1/\alpha)}{n}} + \sigma \frac{\log(U)}{\sqrt{2n\log(1/\alpha)}}\right) \leq \alpha.
\end{equation}
\end{theorem}
The above bound also holds in non-i.i.d.\ settings under a certain martingale dependence structure, where each $X_i$ is $\sigma$-subGaussian conditional on $X_1,\dots,X_{i-1}$, in which case the left-hand side $\EE[X]$ is replaced by $\sum_{i=1}^n \EE[X_i | X_1,\dots,X_{i-1}]/n$. We omit the details for simplicity, see~\cite{howard2020time} for details.

The original Hoeffding inequality~\eqref{eq:hoeff2} is recovered by replacing $U$ with $1$. Since $\log U < 0$, this is a strictly tighter bound than the original. Note that $\EE[\log U] = -1$ and $\Var(\log U)=1$, perhaps giving an intuitive idea of the extent of the gain.

The proof is simple: definition~\eqref{eq:subG} implies that $\exp(\lambda \sum_{i=1}^n (X_i - \EE[X]) - \tfrac{\lambda^2}{2} \sigma^2 n )$ has expected value at most one. Since it is also nonnegative, we can apply our randomized Markov's inequality. Setting $\lambda = \sqrt{\frac{2\log(1/\alpha)}{\sigma^2 n}}$, and rearranging terms, yields the above claim. (Note that $\lambda$ cannot be a function of $U$ for the proof to work.)

\begin{remark}
\label{rmk:hoeffding_length}
The bound in \eqref{eq:hoeff-U} can be rewritten by setting $\epsilon = \sqrt{2\log(1/\alpha)/n}$) as: 
\begin{equation}\label{eq:hoeff-U2}
P\left( \bar X_n - \EE[X] \geq \sigma \epsilon + \frac{\sigma \log U}{\epsilon n} \right) \leq \exp(-n\epsilon^2/2),
\end{equation}
which is easier to compare to~\eqref{eq:hoeff}. Also, a $(1-\alpha)$-confidence interval for the mean $\EE[X]$ is
\begin{equation}\label{eq:hoeff-ci}
\bar X_n \pm  \left[ \sigma \sqrt{\frac{2\log(2/\alpha)}{n}} + \sigma \frac{\log(U)}{\sqrt{2n\log(2/\alpha)}}\right].
\end{equation}
The usual half-width of Hoeffding's inequality simply involves the first term above, that is $\sigma \sqrt{\frac{2\log(2/\alpha)}{n}}$. Since $\EE[\log U]=-1$, the \emph{expected} improvement in width is $\sigma/\sqrt{2n\log(2/\alpha)}$, and the expected \emph{relative} improvement in width is the ratio of the aforementioned quantities (latter divided by former), given by
$1/(2\log(2/\alpha))$. This equals $1/(2\log(40)) \approx 0.14$ when $\alpha=0.05$, meaning we expect to get about a 14\% improvement in width over Hoeffding's inequality. Indeed, we observe exactly this factor of 14\% in our simulations later on.
\end{remark}

\begin{remark}
\label{rmk:bentkus}
The 
Chernoff technique, while attractive due to
its simplicity, leads to inefficient
tail inequalities. 
The original
work of \citeauthor{hoeffding1963}~(\citeyear{hoeffding1963}, Theorem 1) provides 
a simple sharpening of the bound~\eqref{eq:hoeff} for bounded 
random variables. A variety of other strictly tighter
tail bounds
have been developed in the literature---see for 
instance~\cite{bentkus2002,bentkus2006}. Each of these
tail bounds can be derived by
applying MI to a deterministic function of $\bar X_n$, and can therefore be randomized using the UMI in the same way that we have done here. We do not pursue this avenue
further in order to keep our exposition
simple. Nevertheless, in Appendix~\ref{app:bentkus}, we derive
a general inequality which
could be used to derive
such sharper randomized tail bounds,
and which also unifies
Theorems~\ref{thm:umi},~\ref{thm:rand-chebyshev} and~\ref{thm:rand-chernoff}.
\end{remark}

We end by pointing out the possibility that the interval~\eqref{eq:hoeff-ci} can be empty. Despite this, it maintains its frequentist coverage as discussed in Section~\ref{subsec:frequency}, but we nevertheless propose a practical fix in Section~\ref{subsec:truncate}.

\subsection{Exchangeable Hoeffding Inequality} 

If the data are exchangeable rather than i.i.d. $\sigma$-subGaussian, then one can verify via Jensen's inequality that $\bar X_n$ is still $\sigma$-subGaussian. In fact, we have the following stronger claim:

\begin{theorem}[Exchangeable Hoeffding Inequality]\label{thm:exch-hoeff-ineq}
    If $X_1,X_2,\dots$ are an exchangeable sequence of $\sigma$-subGaussian random variables, then for any $\epsilon >0$, 
    \begin{equation}\label{eq:EHI}
    P\left( \sup_{m \geq 1} | \bar X_m - \EE X| \geq \epsilon \sigma \right) \leq \exp(-\epsilon^2/2).   
    \end{equation}
    If the random variables are further assumed to be i.i.d., then we have for any $n \geq 1$,
\begin{equation} 
P\left(\sup_{m\geq n} | \bar X_m - \EE X| \geq \frac{\epsilon\sigma}{\sqrt n}\right) 
\leq \exp(-\epsilon^2/2).
\end{equation}
\end{theorem}

Equation~\eqref{eq:EHI} improves~\eqref{eq:hoeff-one} by recovering it either at $m=1$ or when all random variables are identical. The proof mimics that of Theorem~\ref{thm:exch-chebyshev}, except by using $Y_i = \exp(\lambda (X_i -\EE X) -  \lambda^2\sigma^2/2)$, which has expectation at most one.

\begin{remark}
The above bound is not to be confused with the Azuma-Hoeffding inequality which states that if $X_1,X_2,\dots $ are i.i.d.\ (not merely exchangeable) $\sigma$-subGaussian random variables, then for any $\epsilon > 0$, 
    \begin{equation}\label{eq:azuma-hoeff}
    P\left( \sup_{m \leq n}  \sum_{i=1}^m (X_i - \EE X)  \geq \epsilon \sigma \sqrt{n} \right) \leq \exp(-\epsilon^2/2).   
    \end{equation} 
Of course, the above inequality holds under a certain type of martingale dependence structure, that we omit for simplicity. The Azuma-Hoeffding inequality, despite its fame, is itself loose and was uniformly improved by~\cite{howard2020time}.
\end{remark}

\section{A randomized improvement of Ville's inequality}\label{sec:Ville}

There is a fundamental inequality for nonnegative supermartingales called Ville's inequality, whose time-reversed version lead to the exchangeable Markov inequality. In fact, both Ville's and reverse Ville's inequalities are themselves stronger statements than Markov's inequality.
Here we ask the question: is there a randomized improvement of these inequalities? The answer is subtle---both yes and no in some sense---and we return to answer this question in Section~\ref{sec:Ville-rand} after giving a brief introduction to these inequalities below.

\subsection{Ville's inequality
for forward supermartingales}
\label{sec:Ville-intro}

\cite{Ville:1939} proved that if $X_1,X_2,\dots$ form a nonnegative supermartingale (with respect to any filtration) then for any constant $a>0$,
\begin{equation}\label{eq:ville}
P\left( \sup_{t \geq 1} X_t \geq 1/a\right) \leq a \cdot \EE[X_1].
\end{equation}
Markov's inequality, of course, replaces $\sup_{t\geq 1} X_t$ with just $X_1$, and thus~\eqref{eq:ville}
is strictly stronger. This inequality appears for example in the foundational books of game-theoretic probability~\cite{shafer2005probability,shafer2019game} and much of the recent literature on time-uniform concentration inequalities~\citep{howard2020time,howard2021}, in which it is called Ville's inequality. 

The proof is simple. Define the stopping time $\tau:= \inf\{t \geq 1: X_t \geq 1/a\}$, where $\inf\emptyset=\infty$. For any fixed $m$, Markov's inequality implies
\begin{equation}\label{eq:proof-ville}
    P(\tau \leq m) = P(X_{\tau \wedge m} \geq 1/a) \leq a \cdot \EE[X_{\tau \wedge m}] \leq a \cdot \EE[X_1],
\end{equation}
where the second inequality follows by Doob's optional stopping theorem.
Letting $m\to\infty$ and using the bounded convergence theorem yields $P(\tau < \infty) \leq a\cdot \EE[X_1]$, proving~\eqref{eq:ville}.

\subsection{Randomizing   Ville's inequality}\label{sec:Ville-rand}

\citet[Lemma 3]{howard2021} implies that there are actually three equivalent statements of Ville's inequality: if $M=(M_t)_{t\geq 0}$ is a nonnegative supermartingale with respect to a filtration $\bbF=(\calF_t)_{t\geq 0}$, and $\EE[M_0]=1$, then for any $a>0$, the following three statements hold and imply each other:
\begin{subequations}\label{eq:ville-3versions}
\begin{align}
    P(\exists t \geq 0: M_t \geq 1/a) &\leq a. \label{eq:ville-a}\\
    P(M_\tau \geq 1/a) &\leq a \text{ for every $\bbF$-stopping time }\tau.\\
    P(M_T \geq 1/a) &\leq a \text{ for every } \calF_\infty\text{-measurable random time }T .
\end{align}
\end{subequations}
Despite the fact that the above three statements imply each other, it turns out that (only) the second of these inequalities can be randomized to yield the following result.

\begin{theorem}\label{thm:rand-ville}
For any nonnegative supermartingale $M$ with $\EE[M_0] \leq 1$, we have
\begin{equation}\label{eq:rand-ville}
P(M_\tau \geq U/a) \leq a \text{ for every $\bbF$-stopping time }\tau,
\end{equation}
where $U$ is (stochastically larger than) uniform on $[0,1]$ and is independent of $\bbF$ (and thus independent of $M$ and $\tau$).
\end{theorem}
The proof is simple: the optional stopping theorem for nonnegative supermartingales implies that for any stopping time without restriction, $\EE[M_\tau] \leq 1$, to which we apply the UMI~\eqref{eq:UMI}. 

The entire discussion above does not require time to be discrete: all the statements hold in continuous time as well. In continuous time, it is easy to observe that~\eqref{eq:ville-a} cannot be improved by randomization: the inequality actually holds with equality for $\exp(\lambda B_t - \lambda^2 t/2)$, where $\lambda$ is any nonzero constant and $(B_t)_{t\geq 0}$ is a standard Brownian motion
(see for instance~\cite{durrett2019probability}, Exercise 7.5.2).

We now present a corollary of Theorem~\ref{thm:rand-ville} that has direct implications for sequential testing.

\begin{corollary}[Randomized Ville Inequality]\label{cor:rand-ville}
For any nonnegative supermartingale $M$ with $\EE[M_0] \leq 1$, and any $\bbF$-stopping time $\tau$, we have 
\begin{equation}\label{eq:rand-ville-2}
P(\exists t < \tau: M_t \geq 1/a \text{ or } M_\tau \geq U/a) \leq a, 
\end{equation}
where $U$ is (stochastically larger than) uniform on $[0,1]$ and is independent of $\bbF$ (and thus independent of $M$ and $\tau$).
\end{corollary}

Given any stopping time $\tau$,~\eqref{eq:rand-ville-2} is obtained by applying Theorem~\ref{thm:rand-ville} to the stopping time $\tau' := \min\{\tau, \gamma\}$, where $\gamma :=\inf\{t: M_t \geq 1/a\}$. We now discuss the implications of the above randomized Ville inequality.


\subsection{Implications for safe anytime-valid inference}
\label{sec:implications_rand_ville}

Ville's inequality plays a central role in modern sequential statistics~\citep{howard2020time, howard2021}, and in particular within ``game-theoretic statistics and safe, anytime-valid inference'' \citep{ramdas2020admissible,ramdas2022game}. In the latter context, the above randomized variant of Ville's inequality~\eqref{eq:rand-ville} can improve power in a concrete way. 
To understand why, first note that one usually constructs $M$ to be a nonnegative supermartingale under the null hypothesis, such that it increases to infinity under the alternative. 
Then, $\tau:=\inf\{t: M_t \geq 1/\alpha\}$ is a stopping time of special importance, since Ville's inequality implies that we can reject the null if $\tau<\infty$, while controlling the type-I error at level $\alpha$. 

However, this stopping rule is not literally followed  due to its potentially unbounded nature. One may terminate an experiment even if the above stopping time has not been reached: in a simulation, we have bounded computational resources, which means we only really allow for rejection before some maximum time $t_{\max}$, and in real experiments, one may terminate due (for example) budget constraints. 
This motivates the following three-step rule for sequential testing with the randomized Ville's inequality:
\begin{itemize}
    \item Collect data and continuously monitor the test statistic process $M$ (that is guaranteed to be a nonnegative supermartingale under the null).
    \item If $M$ ever crosses $1/a$, stop and reject the null. Else, stop at any $\mathcal F$-stopping time $\tau$.
    \item Draw an independent random variable $U$ that is (stochastically larger than) uniform on $[0,1]$. Reject the null if $M_\tau \geq U/\alpha$.
\end{itemize}
Theorem~\ref{thm:rand-ville} implies that the above rule yields a bona fide level-$\alpha$ sequential test that is valid under continuous monitoring and adaptive stopping. This is clearly more powerful than the usual rule employed in the aforementioned papers (and references cited therein), which only reject the null when $M_\tau \geq 1/\alpha$.

We note in passing that Theorem~\ref{thm:rand-ville} also applies to a larger class of processes than nonnegative supermartingales. These are called ``e-processes'', and play a particularly key role in sequential composite null testing. However, we omit the details for brevity, and refer instead to the aforementioned survey by~\cite{ramdas2022game}.

We end with the following note. Thanks to a duality between sequential tests and sequential estimation using confidence sequences (CSs), the above observations also have implications for constructing CSs. A CS is a time-uniform or anytime-valid generalization of a confidence interval. Formally, a $(1-\alpha)$-CS for a parameter $\theta$ is a sequence $(C_n(\alpha))_{n \geq 1}$ of confidence intervals (one for each sample size $n$) that are valid at arbitrary stopping times, meaning that it satisfies $P(\theta \in C_\tau(\alpha)) \geq 1-\alpha$ for any $\mathcal F$-stopping time $\tau$, or equivalently it satisfies $P(\forall n \geq 1: \theta \in C_n(\alpha)) \geq 1-\alpha$. Since they are usually (or in fact, essentially always, as per~\cite{ramdas2020admissible,waudby2020estimating}) obtained by inverting a family of sequential tests based on Ville's inequality, our randomized Ville's inequality improves the CI at the final stopping time. To clarify, we have the following result for any CS:
\begin{equation}
\text{For any $\mathcal F$-stopping time $\tau$, we have } P(\exists t < \tau: \theta \notin C_\tau(\alpha) \text{ or } \theta \notin C_\tau(\alpha/U)) \leq \alpha,
\end{equation}
where $U$ is a uniform that is independent of the data and thus the stopping time (one may imagine it to be drawn after stopping). To summarize, if we ever stop a sequential experiment in which we were tracking a CS, the very last confidence interval that we report can be at level $\alpha/U$, and this would still have an overall miscoverage of at most $\alpha$. 



\section{Randomizing the time-reversed Ville inequality}\label{sec:rev-ville}

Directly inspired by~\cite{Ville:1939},
the time-reversed Ville inequality
was first proved
by~\cite{doob1940} for reverse martingales, 
and for instance by \cite{lee1990} and 
\cite{christofides1990} for 
reverse submartingales. We also refer to~\cite{manole2021sequential}
for two self-contained proofs of this result.
The statement is given as follows.
\begin{theorem}[Time-reversed Ville inequality]
\label{thm:reverse_ville}
Let $(X_t)_{t=1}^\infty$ be a nonnegative reverse
submartingale with respect to a reverse filtration
$(\calE_t)_{t=1}^\infty$. Then, for any $a > 0$, 
\begin{equation}\label{eq:rev-ville2}
P\left(\sup_{t \geq 1} X_t \geq 1/a \right) \leq a \cdot \EE[X_1].    
\end{equation}
\end{theorem}
We provide a  self-contained proof
in Appendix~\ref{app:reverse},
which we briefly outline here. 
The proof is similar to that of 
Ville's inequality, with changes to account
for the reversed nature of the process.
Given $m \geq 1$, 
define $\tau:= \sup\{1 \leq t \leq m: X_t \geq 1/a\}$, where $\sup\emptyset=-\infty$. 
Markov's inequality implies
\begin{equation}\label{eq:proof-rev-ville}
    P(\tau \geq 1) = P(X_{\tau \vee 1} \geq 1/a) \leq a \cdot \EE[X_{\tau \vee 1}] \leq a \cdot \EE[X_1].
\end{equation}
To prove the last inequality, notice
that  
the process $Y_t = X_{m-t+1}$, $1 \leq t \leq m$,
is a forward
submartingale, and $\eta:= m-\tau+1$
is a stopping time with respect to 
the same filtration as
$(Y_t)_{t=1}^m$, thus
by the optional stopping theorem,
$$\bbE[X_{\tau \vee 1}] = \bbE[Y_{\eta\wedge m}] 
\leq \bbE[Y_m]=\bbE[X_1].$$
Noting that $\{\tau \geq 1\} = \{\sup_{1 \leq t \leq m} X_t \geq 1/a\}$ yields $P\left(\sup_{1\leq t \leq m} X_t \geq 1/a \right) \leq a \cdot \EE[X_1]$. Sending $m\to\infty$ yields our claim.

\subsection{The uniformly-randomized time-reversed Ville inequality}
Given a reverse filtration $(\calE_t)_{t=1}^\infty$, 
we will say that $\tau$ is a reverse stopping
time if it satisfies $\{\tau = t\} \in \calE_t$
for all $t \geq 1$. We then have the following statement.
\begin{theorem}[Uniformly-randomized time-reversed Ville inequality]
\label{thm:rand_reverse_ville}
Let $(X_t)_{t=1}^\infty$ be a nonnegative reverse
submartingale, and $\tau$ a reverse stopping time, 
both with respect to a reverse filtration
$(\calE_t)_{t=1}^\infty$. Let $U \sim \text{Unif}(0,1)$ be independent of $(\calE_t)$. 
Then, for any $a > 0$, 
\begin{equation}\label{eq:rev-ville2}
P\left( X_\tau \geq U/a \right) \leq a \cdot \EE[X_1].  
\end{equation}
\end{theorem} 
Theorem~\ref{thm:rand_reverse_ville} is a 
uniformly-randomized analogue of 
the time-reversed Ville inequality in Theorem~\ref{thm:reverse_ville}.
This last has been
used for instance by~\cite{manole2021sequential}
to derive nonparametric
sequential goodness-of-fit and two-sample
hypothesis tests using sample convex divergences between
probability distributions as test statistics. 
When the validity of these sequential tests is only needed at arbitrary reverse stopping times $\tau$,
their power can immediately be improved with uniform-randomization by Theorem~\ref{thm:rand_reverse_ville},
in much the same way as we described in 
Section~\ref{sec:implications_rand_ville}.

The proof of Theorem~\ref{thm:rand_reverse_ville}
is straightforward: by reasoning as in the proof
of Theorem~\ref{thm:reverse_ville}, 
above we have $\bbE [X_\tau] \leq \bbE [X_1]$,
thus the claim follows by applying the UMI to $X_\tau$. 

\subsection{A uniformly randomized variant of the EMI}

We can use Theorem~\ref{thm:rand_reverse_ville} to obtain
the following variant of the EMI
which we refer to as the {\it exchangeable
and uniformly-randomized 
Markov inequality (EUMI)}. 
\begin{theorem}[EUMI]\label{thm:umi+emi}
Let $X_1, \dots, X_n$ 
be a set of exchangeable random variables. Then, for any $a \in (0,1)$, 
\begin{align}\label{eq:emi_umi}
    P\left( X_1 \geq U/a \text{ or } \exists t \leq n:  \left|\sum_{i=1}^t X_i/t\right| \geq 1/a \right) \leq a \cdot \bbE |X_1|,
\end{align}
where $U$ is a uniform random variable on $[0,1]$ that is independent of $X_1,\dots,X_n$. 
\end{theorem}
The proof follows by defining the stopping time $\tau := 1\vee  \sup\{1 \leq t \leq n: \sum_{i=1}^t X_i/t \geq 1/a\}$. Then, Theorem~\ref{thm:rand_reverse_ville} implies that
\begin{align}
P\left( \left|\frac 1 \tau\sum_{i=1}^\tau X_i \right|\geq U/a\right) \leq a \cdot \bbE |X_1|,
\end{align}
which is mathematically equivalent to~\eqref{eq:emi_umi}. 
It is easy to see that an analogous statement
holds for infinite sequences of exchangeable
random variables, by removing
the upper bound on $t$ in the definition of $\tau$.  

Equation~\eqref{eq:emi_umi} is one way of
combining the strengths of the EMI and UMI.
Note that the EUMI is stronger than both the UMI and the EMI, and indeed implies both of them.
This inequality has important implications for (more powerful) statistical testing, which we discuss in the following sections. 


\section{Constructing randomized tests with e-values}\label{sec:e-val}

Despite e-values not being defined formally yet, we already used them implicitly in the proofs of all preceding theorems: indeed the random variable $Y_i$ in the proofs of Theorem~\ref{thm:exch-chebyshev} and Theorem~\ref{thm:exch-hoeff-ineq} are e-values, as is $X/\EE[X]$ in the proof of Markov's inequality. 

Before defining e-values formally, we give a brief summary of what's to come, in order to orient the reader. Let $E_1$ and $E_2$ be arbitrarily dependent e-values for testing a given hypothesis. In order to achieve a level $\alpha$ test, the natural way to combine them  is to average them into a combined e-value and check if $(E_1+E_2)/2 \geq 1/\alpha$, whose validity is guaranteed by Markov's inequality. 

However, a simple and uniform improvement is as follows: choose $i \in \{1,2\}$ randomly with equal probability, and first check if $E_i \geq 1/\alpha$, and if not then check whether $(E_1+E_2)/2 \geq 1/\alpha$. This is also a level-$\alpha$ test, due to the EMI. Another valid randomized level-$\alpha$ test is to reject if $(E_1+E_2)/2 \geq U/\alpha$ where $U$ is an independent uniform random variable, whose validity is guaranteed by the UMI. We expand on some statistical applications of these ideas below.

\subsection{A brief review of e-values}

For a set of distributions $\cP$, an e-value\footnote{If $\cP = \{P\}$ is a singleton, then all optimal e-values take the form of $dQ/dP$, that is likelihood ratios of $Q$ against $P$, for some (implicit or explicit) alternative $Q$. So, technically e-values have been around for 100 years masquerading as likelihood ratios (and Bayes factors). The recent christening of the term ``e-value'' is simply to recognize the importance of a more general concept that has utility much beyond point nulls. Indeed, beyond the singleton case, e-values can be viewed as nonparametric/composite generalizations of likelihood ratios to complex settings involving nonparametric and composite nulls and alternatives. Even for this setting, e-values have technically been around for over 50 years~\citep{robbins1970statistical}, appearing in the form of stopped nonnegative supermartingales and implicitly within the proofs of Chernoff bounds~\citep{howard2020time}. The concept appears to have simply floated around without a unified name for 50 years, until several authors---who had a priori used different terms for the same concept---simultaneously decided to converge to the terminology ``e-value'' a few years ago~\citep{Shafer:2021,vovk2021values,GrunwaldHK19,ramdas2020admissible,wasserman2020universal}. Research on e-values has blossomed recently, without acknowledgment of understanding of its roots. The reader may see~\cite{ramdas2022game} for a recent survey on game-theoretic statistics and safe anytime-valid inference, which provides a broader context in which e-values arise naturally, and for many examples of composite, nonparametric e-values, as well as details of the connection to betting scores and the wealth of a gambler betting against the null.}
 for the null hypothesis $H_0: P \in \cP$ is a nonnegative random variable $X$ such that 
\[
\EE_P[X] \leq 1 \text{ for all }P \in \cP.
\]

Since e-values are likely to be small under the null (and hopefully large under the alternative), a level-$\alpha$ test is given by:
\begin{equation}\label{eq:rej-E}
\text{reject $H_0$ if $X \geq 1/\alpha$. }
\end{equation}
This test controls type-1 error (nonasymptotically and hence uniformly over $\cP$) by Markov's inequality~\eqref{eq:M}. Without knowing any further details about the distribution of $X$, this rule does not appear to be improvable in general.



However, we note that with a little external randomization, the rule~\eqref{eq:rej-E} can be made (usually strictly) more powerful. The UMI~\eqref{eq:UMI} implies that if $U$ is an independent uniform random variable on $[0,1]$, then the rule
\begin{equation}\label{eq:rej-E-rand}
\text{reject $H_0$ if $X \geq U/\alpha$, }
\end{equation}
controls the type-I error at level $\alpha$.
In other words, $\min(U/X,1)$ is a valid p-value; this fact was independently recently pointed out in~\cite{ign2022values}. 

We will see that such uniform randomization can effectively be used to improve on (arguably natural) ways to combine e-values to yield tests, but a new type of test is opened up as a consequence of the EMI~\eqref{eq:ex-M}, and the EUMI~\eqref{eq:emi_umi}.

\subsection{Combining multiple arbitrarily dependent e-values to test $\cP$}

Suppose we have constructed $K$ arbitrarily dependent e-values $X_1,\dots,X_K$ for the same null hypothesis $H_0: P \in \cP$. These may or may not be exchangeable. A natural way to form a test is to define 
\[
\bar X_K := (X_1+\dots+X_K)/K,
\]
which is also an e-value, and thus to 
\begin{equation}\label{eq:naive-e-avg}
\text{reject $H_0$ if $\bar X_K \geq 1/\alpha$.}
\end{equation}
Said differently, 
\begin{equation}\label{eq:ep1}
p = \min(1/\bar X_K,1)
\end{equation}
is a p-value.
In fact,~\citet[Appendix G]{vovk2021values} prove that among symmetric e-to-p merging functions, $\min(1/\bar X_K,1)$ is the optimal choice. However, the above optimality result precludes the use of randomization, meaning that it can potentially be dominated by randomized rules. We show that this is indeed the case. To prepare for the result, given a permutation $\pi$ of $\{1,\dots,K\}$, define
\[
\bar X^\pi_t := \frac1t \sum_{i=1}^t X_{\pi(i)}.
\]
Note that $\bar X^\pi_K = \bar X_K$.


 
 \begin{proposition}
     \label{prop:e_value_dependent}
     Let $X_1,\dots, X_K$ be arbitrarily dependent
     e-values for the null hypothesis $H_0:P\in \calP$. If the $X_i$ are not exchangeable, let $\pi$ be a uniformly random permutation of $\{1,\dots,K\}$, otherwise let $\pi$ be the identity permutation.
     Then, the following rules control the type-I error at level $\alpha$: 
     \begin{align}
     \label{eq:u-e-avg}     
     1.& ~ \text{reject } H_0 \text{ if } 
     \bar X_K \geq U/\alpha,  \\
     \label{eq:eval-perm}
     2.& ~ \text{reject } H_0 \text{ if }   \sup_{t \leq K} \bar X^\pi_t \geq 1/\alpha, \\
     \label{eq:u-e-perm-avg}  
     3.& ~ \text{reject } H_0 \text{ if } X_{\pi(1)} \geq U/\alpha \text{ or } \sup_{t \leq K} \bar X^\pi_t \geq 1/\alpha  . 
     \end{align}
     Further, each of these rules is more powerful than~\eqref{eq:naive-e-avg}. 
 \end{proposition}
The proof follows directly from the UMI (Theorem~\ref{thm:umi}) for~\eqref{eq:u-e-avg},
from the EMI (Theorem~\ref{thm:emi}) 
for~\eqref{eq:eval-perm}, and from EUMI (Theorem~\ref{thm:umi+emi}) for~\eqref{eq:u-e-perm-avg}.
Despite
rule~\eqref{eq:eval-perm} being
strictly 
less powerful than rule~\eqref{eq:u-e-perm-avg},
we state it separately for ease of reference below. 
These three rules can be alternatively written as forming one of the following three p-values:
 \begin{align}
 \label{eq:new_ep1}
 1.& ~ p_1 = \min(U/\bar X_K,1) \\ 
 \label{eq:new_ep2}
 2.& ~ p_2 = \min(\inf_{1 \leq t \leq K} (\bar X_t^\pi)^{-1}, 1) \\ 
 \label{eq:new_ep3}
3.& ~ p_3 = \min(\inf_{1 \leq t \leq K} (\bar X_t^\pi)^{-1}, U/X_{\pi(1)}, 1),
\end{align}
and rejecting the null when that p-value is smaller than $\alpha$.

Note that if one wants to combine $K$ arbitrarily dependent e-values to obtain a combined e-value, averaging is still optimal. But if one wants to form a p-value or test at level $\alpha$, options like~\eqref{eq:new_ep1}, \eqref{eq:new_ep2},
\eqref{eq:new_ep3} dominate~\eqref{eq:ep1}.

\subsection{Combining $m$-way independent e-values}

We now briefly extend the previous subsection's ideas beyond arbitrary dependence. Suppose that $X_1,\dots,X_K$ are $m$-way independent e-values for some null $H_0$, meaning that every subset of $m$ of them is jointly independent. $m=1$ corresponds to arbitrary dependence as discussed previously, and $m=K$ corresponds to full joint independence. 

A natural way to test using $m$-way independent\footnote{Even under full independence, such combinations may be more robust and stable than the product $\prod_{i=1}^K X_i$. The latter will have the largest value if all e-values exceed one, but equals zero if even a single e-value was unluckily equal to zero. Of course there are ways to get around the zero issue, like calculating $\prod_{i=1}^K (X_i/2 + 1/2)$, but that is besides the current point.} e-values is to use U-statistics and 
\[
\text{reject if } \bar X_K := \frac1{\binom{K}{m}}\sum_{\substack{A \subset \{1,\dots,K\} \\ \text{s.t. } |A|=m}} \prod_{i \in A} X_i \geq 1/\alpha.
\]
However, a better idea is to randomly sample (with or without replacement) a size-$m$ subset $A_b \subset \{1,\dots,K\}$ one-by-one for $b=1,2,\dots$, calculate the average for the subsets sampled thus far, and reject as soon as any average exceeds $1/\alpha$:
\[
\text{reject as soon as } \frac1{b}\sum_{s=1}^b \prod_{i \in A_s} X_i \geq 1/\alpha \text{ for some $b \geq 1$.}
\]
This is a level-$\alpha$ test because of the EMI~\eqref{eq:ex-M}. Instead, we could also
\[
\text{reject if } \bar X_K := \frac1{\binom{K}{m}}\sum_{\substack{A \subset \{1,\dots,K\} \\ \text{s.t. } |A|=m}} \prod_{i \in A} X_i \geq U/\alpha,
\]
for an independent uniform $U$. As before, while the latter two tests (often strictly) dominate the first one, it is a priori unclear which of the two options will be more powerful, so we explore this in the simulations later. Last, one can use the EUMI to combine the strengths of both the earlier rules, but we omit this for brevity.

\section{Improving the Power of Universal Inference}
\label{sec:univ}

Universal inference~\citep{wasserman2020universal} is a simple and extremely broadly applicable test for any composite null hypothesis, that is nonasymptotically valid without regularity conditions. At its heart is a randomized method for constructing an e-value. To describe the simplest version of their idea, consider a setting where we have i.i.d.\  data $Y_1,\dots,Y_n \sim P$ and we would like to test the null $H_0: P \in \cP$, perhaps against an alternative (implicit or explicit) $H_1: P \in \mathcal Q$. We first partition the data at random into two (possibly unequal) datasets $D_0$ and $D_1$. Using $D_1$, we come up with any estimator/guess $\hat Q$ in $\mathcal Q$. The split likelihood ratio is defined as
\begin{equation}\label{eq:split-lrt}
X =\inf_{P \in \cP}  \prod_{i \in D_0}\frac{ d\hat Q}{dP}(Y_i),
\end{equation}
where we assume for simplicity that $Q \ll P$ for any $P,Q \in \mathcal{P} \cup \mathcal Q$.
In other words, it is the likelihood ratio of a particular alternative $\hat Q$ (picked from $D_1$) against the maximum likelihood estimator under the null. \cite{wasserman2020universal}  prove that $X$ is an e-value for $\cP$, and they
\begin{equation}\label{eq:ui-rule}
\text{reject $H_0$ when } X \geq 1/\alpha.
\end{equation}
There are many other variants, for example using profile likelihoods to handle nuisance parameters, smoothed likelihoods to avoid encountering an infinite likelihood, relaxed likelihoods in case calculating the maximum likelihood is infeasible, and so on. Universal inference is named such because it provides a simple and universally applicable baseline method that works for testing any null, without making any regularity assumptions (unlike the generalized likelihood ratio test, whose threshold is often unknown for singular or complex $\cP$).

Of course, the downside is that the method is conservative in general. In parametric settings without nuisance parameters, when the usual likelihood ratio test applies, universal inference is typically loose (in asymptotic efficiency, say) by a small constant factor of about 2 to 4, achieving the right rate in sample size $n$, dimensionality and level $\alpha$, and appropriate notions of signal-to-noise ratio~\citep{dunn2021gaussian}. In nonparametric settings, sometimes no other test exists, so the conservativeness of universal inference remains unclear~\citep{dunn2021universal}.

We mention two ways to gain back some of the conservativeness. The first is to simply replace~\eqref{eq:ui-rule} with the UMI~\eqref{thm:umi} to yield ``uniformly-randomized universal inference'': let $U$ be an independent $U[0,1]$ random variable, then we may
\begin{equation}\label{eq:rand-univ-inf}
\text{reject $H_0$ when } X \geq U/\alpha,
\end{equation}
yielding a strictly more powerful test than universal inference, that still controls type-I error at level $\alpha$.

The above use of uniform randomization is reminiscent of a somewhat similar use in the context of permutation or randomization tests, which are typically conservative by default, but can be made to be exact by the use of external randomization; see also Section~\ref{subsec:external-rand} on avoiding the use of external randomization $U$.


The second way to gain back some of the constant factors is to recall the randomness associated with sample splitting, and derandomize the approach by averaging. One natural subsampling approach proposed by~\cite{wasserman2020universal} just repeats the calculation of $X$ a total of $B$ times, each time on a different independent random split of the data. Call the resulting exchangeable e-values as $X_1,\dots,X_B$. They propose to 
\begin{equation}\label{eq:UI}
\text{reject $H_0$ when } \frac{X_1+\dots+X_B}{B} \geq 1/\alpha.
\end{equation}
\cite{dunn2021gaussian} prove that such derandomization does improve power.
Given the discussion in preceding sections, using the UMI~\eqref{eq:UMI} to
\begin{equation}\label{eq:UI-U}
\text{reject $H_0$ when } \frac{X_1+\dots+X_B}{B} \geq U/\alpha,
\end{equation}
also yields a level-$\alpha$ test. Alternatively, one may
\begin{equation}\label{eq:UI-M}
\text{reject $H_0$ when } \sup_{1 \leq b \leq B} \frac{X_1+\dots+X_b}{b} \geq 1/\alpha,
\end{equation}
whose validity follows from the EMI~\eqref{eq:ex-M}.

Several of the above observations also apply to other related tests, for example using the reverse information projection e-value~\citep{GrunwaldHK19}.
We summarize the above observations below for easier reference.

\begin{proposition}[(De)randomized universal inference]
Let $X$ be the split likelihood ratio statistic defined in~\eqref{eq:split-lrt} (or the crossfit likelihood ratio defined in~\cite{wasserman2020universal}, or the reverse information projection e-value~\citep{GrunwaldHK19}). Consider the test that
\[
\text{rejects $H_0$ when } X \geq U/\alpha.
\]
Similarly, for subsampling-based universal inference, consider the rules~\eqref{eq:UI-U} or~\eqref{eq:UI-M}, or a third rule that
\begin{equation} 
\label{eq:eumi_ui}
\text{rejects $H_0$ when either } X_1 \geq U/\alpha \text{ or } \sup_{1 \leq b \leq B} \frac{X_1+\dots+X_b}{b} \geq 1/\alpha.
\end{equation}
All the above rules are  more powerful than universal inference~\eqref{eq:ui-rule}, and control type-I error at level $\alpha$ without regularity conditions.
\end{proposition}

An important point worth remarking is that $B$ has to be fixed in advance if we want to average the e-values to calculate a single e-value. However, in~\eqref{eq:UI-M} or~\eqref{eq:eumi_ui}, $B$ \emph{does not have to be fixed in advance} since technically those suprema hold from 1 to $\infty$, 
not just 1 to $B$ as stated, and thus they hold even if $B$ is chosen adaptively (as a stopping time, say). In other words, instead of fixing $B$ in advance, one can simply calculate one $X_b$ at a time, calculate running averages as you go along, stop when you want, and reject if at any step the average crosses $1/\alpha$ (or, in the case
of~\eqref{eq:eumi_ui}, stop as early as the first
step if $X_1 \geq U/\alpha$).


\section{(De)randomizing nonparametric tests based on betting}
\label{sec:betting}
Consider a simple special case of a nonparametric testing problem from~\cite{waudby2020estimating}. Let $Y_1,\dots,Y_n$ be drawn i.i.d.\  from an unknown distribution $P$ on $[0,1]$, having mean $\mu$. Suppose we want to test
the null $H_0: \mu = 0.5$, against an alternative $H_1: \mu > 0.5$. Define the initial wealth of a gambler who wishes to bet against this null as $M_0 = 1$, and let their wealth evolve as
\begin{equation}
\label{eq:Mt}
M_t = \prod_{i=1}^t (1+ \lambda_i (Y_i - 0.5)) = M_{t-1}\cdot (1+\lambda_t(Y_t-0.5)),
\end{equation}
where $\lambda_i \in [0,2]$ is a random variable (representing the gambler's bet) that can be chosen based on $Y_1,\dots,Y_{i-1}$, meaning that it is ``predictable'' with respect to the filtration $\mathcal F_t:=\sigma(Y_1,\dots,Y_t)$. It is easy to check that under the null, $(M_t)_{t \geq 0}$ is a nonnegative martingale with initial value one (and in fact for each fixed $t$, $M_t$ is an e-value). 
Ville's inequality 
(Section~\ref{sec:Ville-intro}) 
implies that, under the null,
$P\left( \sup_{0 \leq t \leq n} M_t \geq 1/\alpha \right) \leq \alpha,$
and thus a level-$\alpha$ test is obtained by
\begin{equation}\label{eq:bet-ville2}
\text{rejecting the null if } \sup_{0 \leq t \leq n} M_t \geq 1/\alpha.
\end{equation}
Said differently, $\inf_{0 \leq t \leq n} (1/M_t)$ is a p-value. The experiments in \cite{waudby2020estimating} demonstrate that these tests perform excellently in practice (when using appropriate rules\footnote{One can show that if the alternative is true, then it is possible to bet smartly (meaning derive automated rules to predictably set $\lambda_t$) so that the gambler's wealth $M_t$ grows exponentially fast, with the exponent automatically adapting to both the unknown signal $\mu-0.5$, and the unknown variance $\EE[(Y-\mu)^2]$. The authors also derive new exponential ``empirical Bernstein'' inequalities that can achieve the same effect.} for updating $\lambda_i$ at each step), and the resulting confidence intervals obtained by inverting such a test are usually much shorter than a plethora of competing methods. We note immediately an improvement delivered by the randomized Ville inequality: we may
\begin{equation}\label{eq:bet-ville3}
\text{reject the null if } \sup_{0 \leq t \leq n} M_t \geq 1/\alpha \text{ or } M_n \geq U/\alpha,
\end{equation}
for an independent uniform random variable $U$.

However, there is something slightly unsettling about this test: 
the p-value and test
depend on the random order $Y_1,\dots,Y_n$ of processing the points one by one, since the bets $\lambda_t$ depend on $Y_1,\dots,Y_{i-1}$. Since the chosen ordering was random to begin with, the p-value and test are symmetric functions of the data in a distributional sense (indeed, one can randomly scramble the data before running the test to enforce this), but there is some sense in which one may hope that the ``algorithmic randomness'' introduced by processing the data along \emph{one} random ordering method can somehow be removed. In fact, this is also an issue with the oracle test: there is an optimal choice $\lambda^*(P)$ that could be used at every step, but the supremum over $t$ still makes the resulting test or p-value dependent on the order of processing the bag of data.

\cite{waudby2020estimating} describe one way to remove the effect of the arbitrary ordering. Noting that the final wealth $M_n$ is an e-value (because of the nonnegative martingale property of the wealth process under the null), one can repeat the procedure $B$ times on different independent random permutations of the original data, and only note the final wealth at time $n$ on each permutation of the data, denoted $M^1_n,\dots,M^B_n$.
Then, they proposed to 
\begin{equation} 
\label{eq:betting_MI}
\text{reject the null if } \frac{M^1_n + \dots + M^B_n}{B} \geq 1/\alpha,
\end{equation}
which for large enough $B$ effectively becomes a symmetric function of the data.

Noting the parallels between the above rule and~\eqref{eq:UI}, one can instead gain more power by using~\eqref{eq:UI-U} or~\eqref{eq:UI-M} instead, meaning to either
\begin{equation}\label{eq:bet-U}
\text{reject the null if } \frac{M^1_n + \dots + M^B_n}{B} \geq U/\alpha,
\end{equation}
or to
\begin{equation}\label{eq:bet-M}
\text{reject the null if } \sup_{1\leq b\leq B}\frac{M^1_n + \dots + M^b_n}{b} \geq 1/\alpha.
\end{equation}
The first method derandomizes by averaging (thus removing the effect of data ordering), but then again randomizes the threshold using the UMI. The second derandomizes in a more sophisticated manner, using the EMI. The following rule combines the two techniques
\begin{equation}\label{eq:bet-umi-emi}
\text{reject the null if } M_n^1 \geq U/\alpha \text{ or } \sup_{1\leq b\leq B}\frac{M^1_n + \dots + M^b_n}{b} \geq 1/\alpha,
\end{equation}
using the EUMI.
We formalize these observations below.
\begin{proposition}
\label{prop:betting}
Let $Y_1, \dots, Y_n$ be 
independent random
variables supported in [0,1], with identical mean $\mu$. 
Let $\pi_1, \dots, \pi_B$ be permutations of 
$\{1, \dots, n\}$ that are sampled uniformly at random (with or without replacement), and let
$$M_n^b = \prod_{i=1}^t (1+\lambda_i(Y_{\pi_b(i)}-0.5)),
\quad b=1, \dots, B.$$
Then, the test~\eqref{eq:bet-ville3} 
for $H_0: \mu = 1/2$ controls the type-I error
at level $\alpha$, and is more powerful than the  original rule~\eqref{eq:bet-ville2}.
Likewise, the tests~\eqref{eq:bet-U},~\eqref{eq:bet-M} and~\eqref{eq:bet-umi-emi}
control the type-I error
and are more powerful than the rule~\eqref{eq:betting_MI}. 
\end{proposition}

Since one cannot take a supremum over both $t$ and $B$, it is a priori unclear which of the tests~\eqref{eq:bet-ville2} and \eqref{eq:bet-M} is more powerful (or~\eqref{eq:bet-ville3} versus~\eqref{eq:bet-umi-emi}). We examine such questions in the simulations that follow.

We note that when inverting these tests to form confidence intervals for the mean $\mu$, as done in~\cite{waudby2020estimating}, the same $U$ can be used across all the tests.

\section{Experiments}
\label{sec:simulations}
We perform a simulation
study to illustrate the extent to which
the UMI, EMI, and EUMI can increase the power of the
aforementioned methodologies. 
Code for reproducing this simulation
study is publicly
available\footnote{\href{https://github.com/tmanole/Randomized-Markov}{https://github.com/tmanole/Randomized-Markov}}. We choose the level $\alpha=.05$
across all simulations. The parameter~$B$ appearing
in Sections~\ref{sec:univ}--\ref{sec:betting}
is always taken to be 100.


\subsection{Confidence Intervals for a
Gaussian Mean}
\begin{figure}[htbp!]
\centering
\begin{minipage}{0.45\textwidth}
\includegraphics[width=\textwidth]{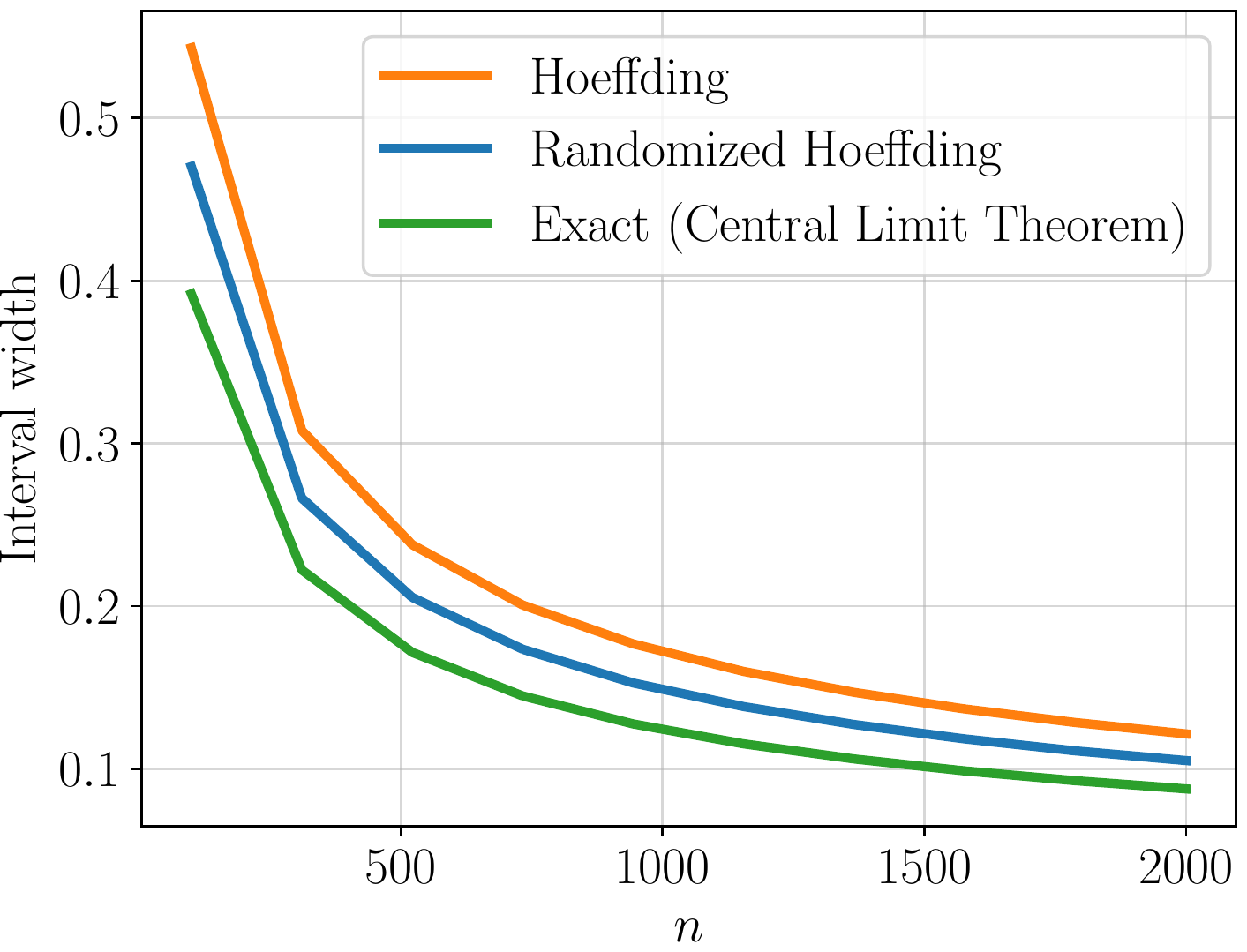}
\caption*{ (a) Average confidence interval length.}
\end{minipage}~~~~~
\begin{minipage}{0.45\textwidth}
\includegraphics[width=1.04\textwidth]{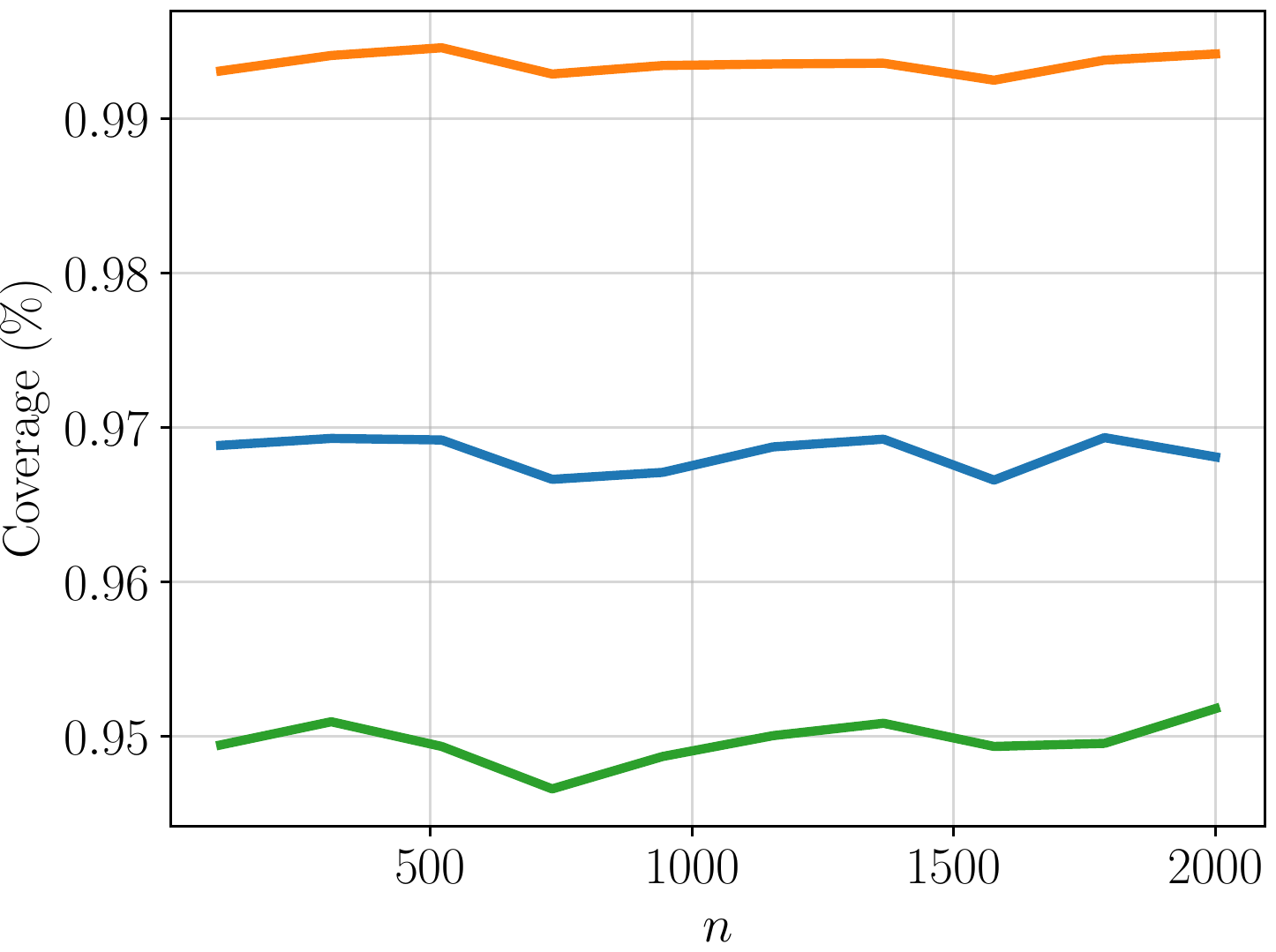}
\caption*{(b) Average confidence interval coverage.}
\end{minipage}
\caption{\label{fig:ci_hoeff} Average length and coverage 
of the three
confidence intervals across 
20,\,000 replications for each sample size.
The uniformly-randomized Hoeffding interval sits halfway
between the traditional Hoeffding interval and the exact
interval, both in terms of coverage and length.}
\end{figure}

We begin with a toy example to 
compare the 
tightness of our randomized tail bounds. 
Let $X_1, \dots, X_n$ be an i.i.d.
sample from the $\calN(0,1)$
distribution. We
compare the width of the 
uniformly-randomized Hoeffding 
confidence interval~\eqref{eq:hoeff-U} 
for the mean $\bbE[X_1]$, to that
of the traditional Hoeffding 
interval~\eqref{eq:hoeff2}. 
As a benchmark, we also compare
them to  the exact confidence interval
$\bar X_n \pm z_{\alpha/2}/\sqrt n$,
where $z_{\alpha/2}$ is the $1-\alpha/2$
quantile of the standard Gaussian distribution.
The average length and coverage of these
three intervals is reported in Figure~\ref{fig:ci_hoeff},
across ten values of $n \in [100,\!2000]$.

It can be seen that the 
randomized Hoeffding
interval (based on the UMI)
has length lying between
that of the traditional Hoeffding
interval and the exact interval. 
By reasoning similarly as in Remark~\ref{rmk:hoeffding_length}, 
the expected relative
improvement in length
of the exact interval over the 
UMI interval is approximately 17\%.
On the other hand, we have already
stated that the relative
improvement in length of the 
randomized Hoeffding
bound over its
classical counterpart is 14\%. Both of these 
expected length
ratios are confirmed
by our simulation study,
and show that the UMI interval sits roughly
halfway between the Hoeffding and exact
intervals, both in length and 
coverage. 
As discussed in Remark~\ref{rmk:bentkus},
the Hoeffding
interval can itself 
be sharpened using more sophisticated
tail bounds, and we expect
that randomized versions of 
such inequalities
would lead to even tighter intervals.

\subsection{Testing with a set of arbitrarily dependent e-values}
\label{sec:sim_evals}

Assume that $X_1,\dots,X_K$ are $K=100$
test statistics which are not necessarily
independent, and which are distributed
as $ \calN(\mu,1)$ for some $\mu \in \bbR$.
We would like to combine them
to test the null hypothesis $H_0:\mu \leq 0$. One approach is to define the transformed statistics
\begin{equation} 
\label{eq:sim-e-vals}
E_j = \exp(X_j - 0.5), \quad j=1,\dots,K,
\end{equation}
which are e-values for the null hypothesis $H_0$. 
We can
combine them to test $H_0$ using 
the rule~\eqref{eq:naive-e-avg}
based on averaging the e-values and applying MI (``Av+MI''), the rule 
\eqref{eq:u-e-avg} based on  UMI, 
the rule~\eqref{eq:eval-perm}
based on  EMI, 
or the rule~\eqref{eq:u-e-perm-avg} based on
on EUMI.

\begin{figure}[h!]
\centering
\includegraphics[width=0.9\textwidth]{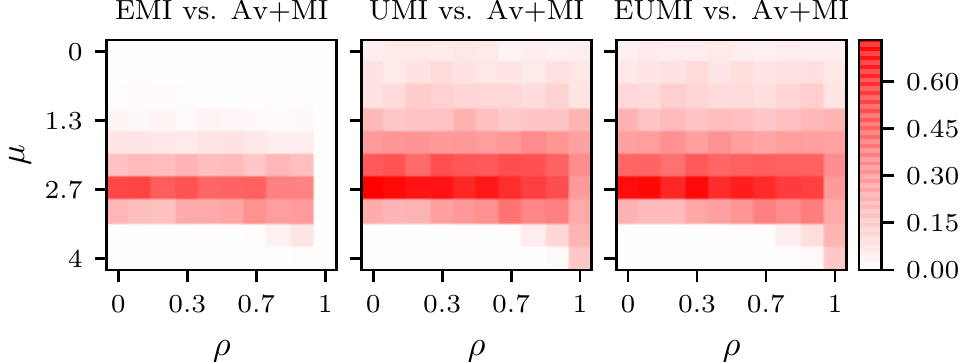}
\caption{\label{fig:e-val} Comparison of the
rejection proportions $\pi_{\mathrm{AvMI}}$, $\pi_{\mathrm{EMI}}$, $\pi_{\mathrm{UMI}}$,
and $\pi_{\mathrm{EUMI}}$. 
of the procedures~\eqref{eq:naive-e-avg}, \eqref{eq:u-e-avg}, ~\eqref{eq:eval-perm},
and~\eqref{eq:u-e-perm-avg}, 
for rejecting the null hypothesis~$H_0:\mu\leq 0$
based on the e-values~\eqref{eq:sim-e-vals}.
For varying
values of~$\rho$ and~$\mu$, the left-hand
side plot represents 
the difference $\pi_{\mathrm{EMI}} - \pi_{\mathrm{AvMI}}$, 
the middle plot represents
$\pi_{\mathrm{UMI}} - \pi_{\mathrm{AvMI}}$, 
and the right-hand side plot represents 
$\pi_{\mathrm{EUMI}} - \pi_{\mathrm{AvMI}}$.
The procedures based on the UMI, EMI, and EUMI
are strictly more powerful
than the naive procedure based on MI, in some
cases leading to an absolute increase in power of 60\%.}
\end{figure}

We simulate the power of these four
approaches by drawing $(X_1, \dots, X_K)$ with a Toeplitz-structured covariance matrix, where
$\text{\bf Cov}[X_i,X_j] = \rho^{|i-j|}$
for all $1 \leq i,j\leq K$ and some $\rho \in [0,1]$. 
In Figure~\ref{fig:e-val}, we 
report their proportion of rejections 
across ten equally-spaced 
values of the true mean $\mu \in [0,4]$  
and   of the correlation parameter $\rho \in [0,1]$. 
The results are based on 500 replications
from each model under consideration.

It can be seen that the tests based on the EMI, UMI,
and EUMI
improve upon the test based on the MI by at least 10\% (in absolute power)
for most combinations $(\mu,\rho)$. In some
cases, the improvement in power is as high as 60\%
for the UMI and EUMI. The gain made by the UMI 
is similar to that of the EUMI, and both
have similar behavior across all values of $\rho$.
In contrast, the gains made by the EMI are most pronounced for small values of $\rho$.
This is to be expected, since the cumulative averages of e-values in equation~\eqref{eq:eval-perm} are highly correlated
when $\rho$ is large; in fact, they are identical
in the limit $\rho=1$, in which case the 
rejection rules based on MI and EMI coincide.

We  report the individual rejection
proportions of all four methods
in Figure~\ref{fig:e_val_all}
of Appendix~\ref{app:more_sims}.
Therein, we also
report simulation results 
under $K=2$ rather than $K=100$; see Figures~\ref{fig:K2-e-val}--\ref{fig:K2-e_val_all}.
When $K=2$, it can be seen that
the EUMI provides a more pronounced improvement
over the UMI, with an absolute gain in power
as high as 10\% for several combinations $(\mu,\rho)$.

\subsection{Universal Inference for model selection}
We compare the methods presented
in Section~\ref{sec:univ} for the problem
of testing the number
of components in a Gaussian mixture model. 
It is well-known that 
the parametric
family of Gaussian mixtures does not
satisfy the regularity conditions required
for (twice the negative logarithm of) the likelihood ratio statistic to admit
its traditional $\chi^2$ limiting distribution~\citep{ghosh1984asymptotic,dacunha1999testing,chen2009hypothesis}. 
In contrast, the method of universal inference
based on the {\it split} likelihood ratio statistic, 
and its variants presented in Section~\ref{sec:univ},
are valid without any regularity conditions, 
and are therefore natural candidates for this problem. 

Although the limiting distribution 
of the likelihood ratio statistic is unknown
or intractable for general Gaussian mixtures, 
it admits a simple expression when the underlying
mixing 
proportions are known~\citep{goffinet1992testing}. 
We will assume this to be the case so that
we can use the likelihood ratio test (LRT) as a benchmark, 
but we emphasize that the LRT cannot easily
be used to derive a valid test for more
general Gaussian mixtures, where the universal inference method would remain valid. 

Let
\begin{equation} 
\label{eq:mixture} 
X_1, \dots, X_n \overset{\mathrm{i.i.d.}}{\sim} 
0.25\cdot \calN(\mu_1,1) + 0.75\cdot \calN(\mu_2,1),
\end{equation}
where the only
unknown parameters are $\mu_1,\mu_2 \in \bbR$, and consider the problem of testing
whether the above mixture has one vs. two components, i.e.
\begin{equation} 
\label{eq:hypothesis_mixture} 
H_0: \mu_1=\mu_2, \quad \text{vs.}\quad H_1: \mu_1 \neq \mu_2.
\end{equation}
By Theorem~1 of \cite{goffinet1992testing}, if $\lambda$
denotes the likelihood-ratio statistic
for these hypotheses, then $-2\log \lambda$ 
admits the limiting distribution $\max(0,Z)^2$,
for $Z \sim \calN(0,1)$, thus a valid level-$\alpha$
test for $H_0$ is to 
$\text{reject if } -2\log \lambda > q_{1-2\alpha},$
the latter being the $1-2\alpha$ quantile
of the $\chi_1^2$ distribution. 
We will refer to this as the LRT test. We compare its numerical
performance to that of
universal inference (UI; \eqref{eq:ui-rule}), 
uniformly-randomized UI (UMI-UI;~\eqref{eq:rand-univ-inf}),
subsampling UI (SUI;~\eqref{eq:UI}), 
uniformly-randomized SUI (UMI-SUI;~\eqref{eq:UI-U}), 
exchangeable SUI (EMI-SUI;~\eqref{eq:UI-M}),
and exchangeable, uniformly-randomized
SUI (EUMI-SUI;~\eqref{eq:eumi_ui}). 
\begin{figure}[htbp!]
\centering
\includegraphics[width=0.5\textwidth]{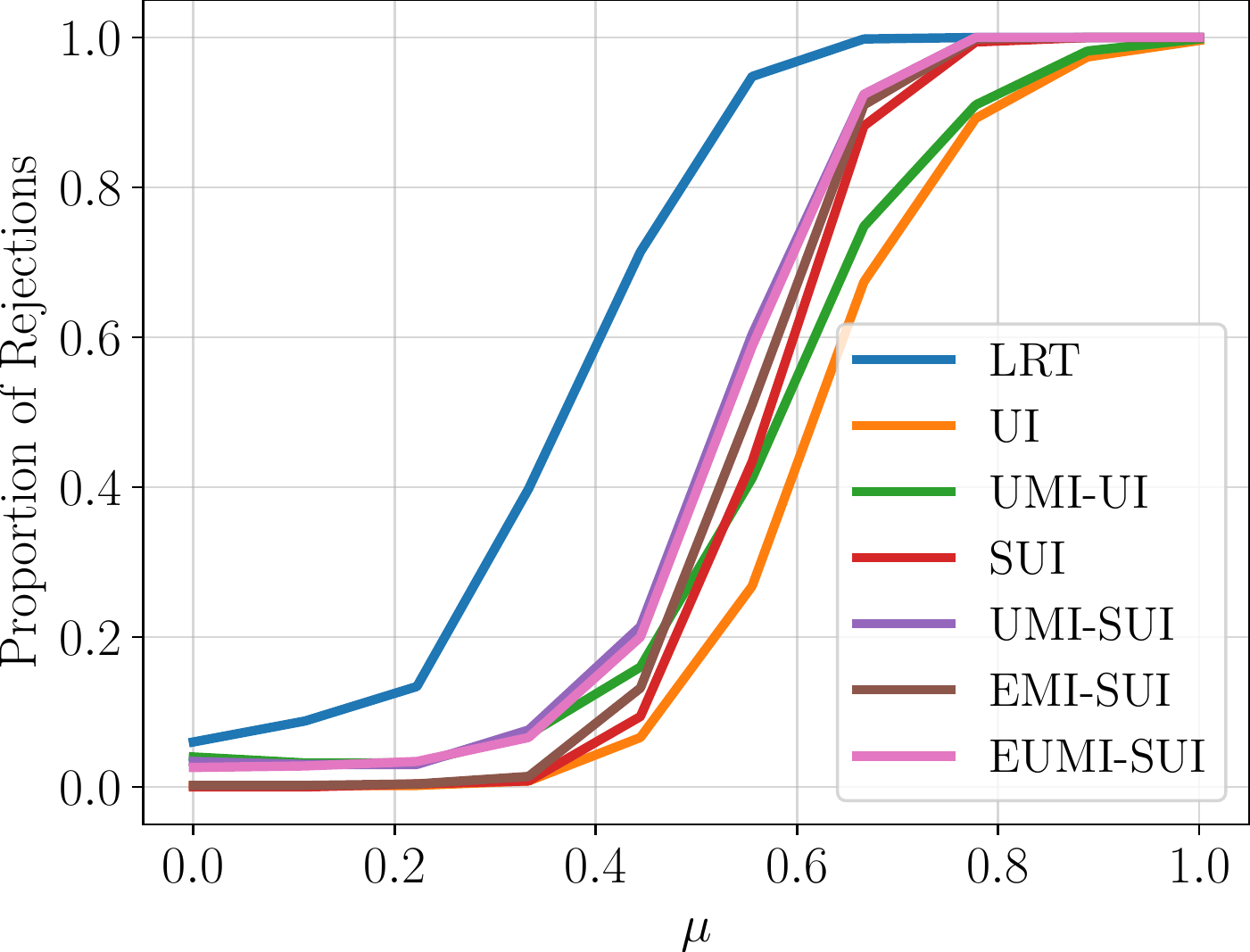}
\caption{\label{fig:ui} 
Empirical power of the seven tests
for the null hypothesis $H_0$ in equation~\eqref{eq:hypothesis_mixture}.
While the benchmark test provided by the LRT
is most powerful, the variants of Universal Inference
based on the EMI, UMI, and EUMI are more powerful
than their MI counterparts. 
}
\end{figure}

Figure~\ref{fig:ui} 
reports the empirical power
of these procedures based on 500 samples
of size $n=500$ from model~\eqref{eq:mixture}, 
under ten equally-spaced values of $\mu := -\mu_1 = \mu_2 \in [0,1]$.  
We observe that the power
of the UMI-UI method
uniformly dominates that of
the original UI method, 
and similarly, the UMI-SUI, EMI-SUI,
and EUMI-SUI methods dominate their SUI counterpart.
The methods UMI-SUI and EUMI-SUI 
exhibit similar performance, and
their absolute increase in power compared
to SUI is on the order
of 15\% for some values of $\mu$. 
Although all methods based on the split
LRT are markedly more conservative
than the LRT, we recall that 
they can be used in arbitrary mixture models,
while the LRT cannot.

\subsection{Testing the mean of a bounded random variable by betting}
\label{sec:sim_betting}

Let $X_1, \dots, X_n$ be an i.i.d. sample from 
a $\text{Beta}(a,b)$ distribution, for some $a,b > 0$. 
Let $\mu = \bbE[X_1] = a/(a+b)$, 
and consider the problem of testing the null hypothesis 
$H_0:\mu=1/2$ using the procedures
defined in Section~\ref{sec:betting}. We form the wealth
statistic $M_n$ in equation~\eqref{eq:Mt} 
based on a 
predictable sequence $(\lambda_i)_{i=1}^n$ chosen according to the LBOW
betting strategy described in~\cite{waudby2020estimating}, 
and compare the rejection 
rules in equations~\eqref{eq:bet-ville2},~\eqref{eq:betting_MI},~\eqref{eq:bet-U},~\eqref{eq:bet-M},\eqref{eq:bet-umi-emi},
based respectively on Ville's
inequality, averaging followed by MI (``Av+MI''), UMI, EMI, and EUMI. 
We take
$a=20$, thus the null hypothesis reduces
to $H_0: b=20$, and we compare the empirical power of these methods 
for varying values of $b \in [19,20.8]$ and $n \in [100,2000]$, based on 500 replications for each pair~$(n,b)$.

\begin{figure}[h!]
\centering
\includegraphics[width=0.78\textwidth]{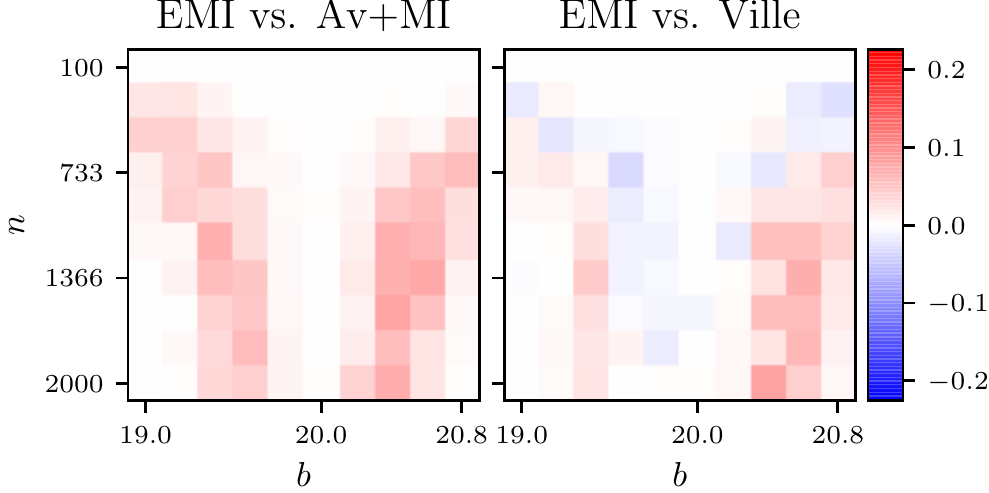}
\caption{\label{fig:betting_emi} Comparison of the
rejection proportions $\pi_{\mathrm{AvMI}}$, $\pi_{\mathrm{Ville}}$, and $\pi_{\mathrm{EMI}}$ 
of the respective procedures~\eqref{eq:betting_MI}, \eqref{eq:bet-ville2} and~\eqref{eq:bet-M},
for rejecting the null hypothesis~$H_0:b=20$
based on the statistic $M_n$.
For varying
values of~$b$ and~$n$, the left-hand
side plot represents 
the difference $\pi_{\mathrm{EMI}} - \pi_{\mathrm{AvMI}}$, 
and the right-hand side plot represents 
the difference $\pi_{\mathrm{EMI}} - \pi_{\mathrm{Ville}}$.
The procedure based on the EMI provides a modest
improvement over
that based on MI, but does not dominate
the procedure based on Ville's inequality.
} 
\end{figure}

\begin{figure}[h!]
\centering
\includegraphics[width=0.78\textwidth]{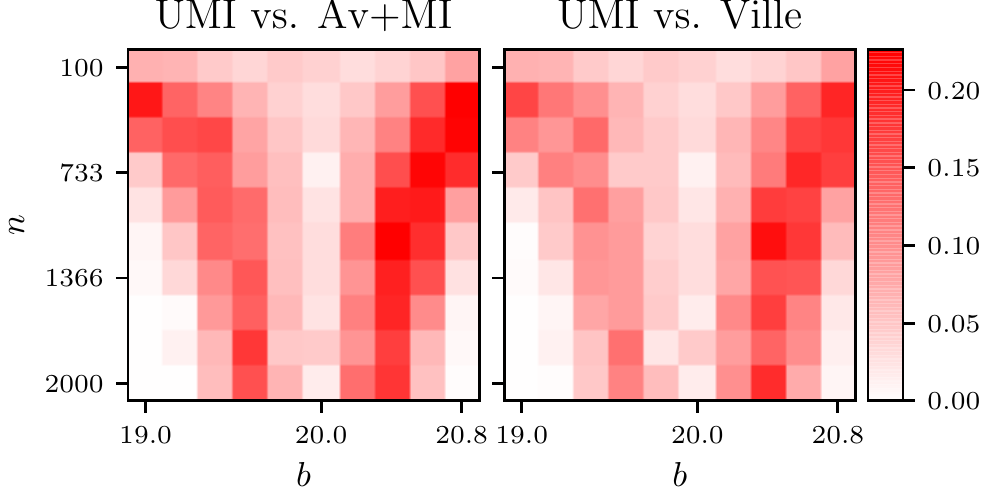}
\caption{\label{fig:betting_umi} Comparison of the
rejection proportions $\pi_{\mathrm{AvMI}}$, $\pi_{\mathrm{Ville}}$, and $\pi_{\mathrm{UMI}}$ 
of the respective procedures~\eqref{eq:betting_MI}, \eqref{eq:bet-ville2} and~\eqref{eq:bet-U},
for rejecting the null hypothesis~$H_0:b=20$
based on the statistic $M_n$.
For varying
values of~$b$ and~$n$, the left-hand
side plot represents 
the difference $\pi_{\mathrm{UMI}} - \pi_{\mathrm{AvMI}}$, 
and the right-hand side plot represents 
the difference $\pi_{\mathrm{UMI}} - \pi_{\mathrm{Ville}}$.
In contrast to the EMI-based procedure reported
in Figure~\ref{fig:betting_emi}, the UMI-based procedure
uniformly dominates both procedures based on MI and Ville's inequality.}
\end{figure}

Figure~\ref{fig:betting_emi} 
compares the EMI-based procedure
to those based on the MI and Ville's inequality. 
It can be seen that the EMI
yields a modest improvement
in power---on the order of
5\%---compared to the 
MI, across the majority
of choices of $b$ and $n$. 
In some cases, it yields
an improvement of similar order over the 
procedure based on Ville's
inequality, but does not 
uniformly dominate
this method, as could have
been anticipated from the discussion
in Section~\ref{sec:betting}. 
In contrast, in Figure~\ref{fig:betting_umi},
it can be seen that the procedure
based on the UMI uniformly dominates
both the MI and the procedure based
on Ville's inequality,
with a gain in absolute power as high as 20\% in 
many cases. The performance
of the EUMI-based procedure is nearly
identical to that of the UMI, thus
we defer this result to Figure~\ref{fig:betting_eumi}
of Appendix~\ref{app:more_sims}.

In Figure~\ref{fig:betting_all} of Appendix~\ref{app:more_sims}, 
we report the individual power
of these five methods. In particular, 
it can be seen therein that the 
power of all methods is at most approximately
5\% under the null hypothesis $b=20$. It can further
be seen that all methods achieve perfect (or nearly perfect) power for values of $b$ near the boundaries
of the interval $[19, 20.8]$, and for large
values of $n$. This explains why 
the UMI, EMI, and EUMI do not to
appear to provide a substantial
improvement over the other methods in this regime.

\section{Discussion}\label{sec:sum}

This paper presented the uniformly-randomized and exchangeable Markov inequalities, along with extensions to Chebyshev and Chernoff bounds, and some example statistical applications involving universal inference and testing by betting. We now begin a relatively long discussion, mixinsg technical observations with some philosophical thoughts.

\subsection{Markov's inequality as a derandomization of UMI}

Let us begin by recalling  the following implication of Markov's inequality from Section~\ref{sec:e-val}. For a nonnegative, integrable $X$, we have that $X/\EE[X]$ is an e-value, and thus $p^* := \EE[X]/X$ is a p-value, meaning that
\[
P(\EE[X]/X \leq a) \leq a.
\]
In this vein, UMI implies that $U \EE[X]/X$ is also a valid p-value. 

What if one tries to derandomize this statement? Suppose we draw $B$ independent uniforms $U_1,\dots,U_B$ and calculate $B$ such p-values (where $p_b = U_b \EE[X]/X$). Then, we get $B$ dependent, exchangeable, p-values. An old result by \cite{ruschendorf1982random} implies that twice the average of arbitrarily dependent p-values are also p-values, meaning that $\bar p_B := 2(p_1 + \dots + p_B)/B$ is a p-value, and this factor of 2 cannot be improved in general (also see~\cite{vovk2020combining}). \cite{choi2023averaging} recently showed that the factor of 2 cannot be improved even assuming that the p-values are exchangeable (which is true in our setting). 

Now note that as $B \to \infty$, $\bar p_B$ converges to $p^*$, because the uniforms average out to 1/2. Thus the gain made by the uniform randomization is exactly offset by the factor of 2 lost by combining p-values. In other words, one can view Markov's inequality as a derandomized version of our UMI.

Further note that other forms of randomization do not appear to help. A result by \cite{ruger1978maximale} shows that $\tilde p_B := 2 \cdot\text{median}(p_1,\dots,p_B)$ is also a p-value and the factor of 2 cannot be improved (assume $B$ is odd for simplicity). In fact, the same work also showed that $p^k_B := p_{(k)} B/k$ is a p-value for any fixed $k \leq B$, where $p_{(k)}$ is the $k$-th smallest p-value. Remarkably, $\tilde p_B$ and $p^k_B$ also converge to $p^*$ as $B \to \infty$. This appears to be a perfectly-designed coincidence, but perhaps on more reflection a simple explanation of this phenomenon may be found. For now, it adds further justification to the title of this subsection.

\subsection{On the role of external randomization in statistics}

At a high level, there appear to be (at least) three reasons that external randomization is used in statistics:
\begin{enumerate}
    \item \textbf{To save computation.} A classic example of this would be the permutation test. When applied to (say) a problem like two-sample or independence testing, the deterministic permutation test needs $n!$ permutations, where $n$ is the number of data points. The variant that is typically used in practice, however, involves the permutations being uniformly sampled from the set of all permutations. This also results in a valid p-value, and the randomization is introduced solely to save computational effort. Another example is a risk-limiting election audit, where ballots are sampled in a random order, allowing one to possibly stop the audit early (with a guarantee on the error); without this randomization, one must look at every single ballot in the audit. A last example would be the use of stochastic gradient descent in minimizing convex objectives because calculating a full gradient may be too expensive.
    \item \textbf{To enable inference that is (essentially) impossible otherwise.} An example of this would be differential privacy; it is provably impossible to guarantee privacy without randomization (adding noise to summary statistics before releasing them). The same impossibility also arises in online learning against adversaries (e.g.: adversarial multi-armed bandits), and in calibrated probabilistic forecasting. Such examples also arise in the Monte Carlo literature, for example MCMC algorithms are used to enable sample from distributions that would otherwise be (analytically and computationally) intractable. The bootstrap and subsampling methods would also be examples. Universal inference also falls under this umbrella, since for many problems, we do not know of any other computational or analytical tool that could replace it. Another example is the knockoffs method for (fixed-X or model-X) conditional independence testing, and the related conditional randomization test. Nonexchangeable conformal prediction provides a contemporary case in point. Finally, sample splitting is used in a variety of contexts to enable assumption-lean inference (such as in post-selection inference).
    
    \item \textbf{For more powerful inference.} An example of this would again be the permutation test. Whether used in its deterministic or random form, it produces a discrete p-value. Randomization can be introduced to convert it into a continuous p-value that is almost surely smaller than the discrete p-value (and exactly uniformly distributed under the null), thus improving (strictly) power.
\end{enumerate}
The randomization used by UMI falls into the last category: the only purpose of introducing randomization is to improve power. 

\bigskip
On a different note, there are at least two types of randomized procedures:
\begin{enumerate}
    \item \textbf{Can be computationally derandomized.} For many procedures, expending more computation can render them ``effectively deterministic'' in the sense that some concentration of measure kicks in, so that the stochastic result concentrates around some deterministic limiting quantity. Examples include the bootstrap, subsampling, permutation and  Monte Carlo methods (including MCMC), the conditional randomization test, (subsampling-based) universal inference, and stochastic gradient descent for convex optimization.
    \item \textbf{Cannot be computationally derandomized.} This includes  adversarial multi-armed bandits, probabilistic forecasting, differentially private inference, nonexchangeable conformal prediction, and many methods based on sample splitting.
\end{enumerate}

UMI falls into the first category, but as discussed in the previous subsection, the power benefits of UMI vanish when it is derandomized, because it reduces exactly to Markov's inequality.

\subsection{Frequency interpretation of the tests and confidence intervals}\label{subsec:frequency}
Chebyshev's inequality~\eqref{eq:cheby} implies that for any $\alpha \in (0,1)$,
\begin{equation}\label{eq:cheby-ci}
  \bar X_n \pm \frac{\sigma}{\sqrt{\alpha n}}   
\end{equation}
is a $(1-\alpha)$ confidence interval for $\EE X$. This means that when we construct infinitely many such intervals for different problems (with independent data from distributions with potentially different means and variances), at least 95\% of those confidence intervals will cover the corresponding means. In contrast, our randomized Chebyshev's inequality~\eqref{eq:cheby} implies that 
\begin{equation}\label{eq:cheby-u-ci}
  \bar X_n \pm \frac{\sigma\sqrt{U}}{\sqrt{\alpha n}}   
\end{equation}
is also $(1-\alpha)$ confidence interval for $\EE X$. We highlight that it has \emph{exactly} the same frequency interpretation as above. Despite being randomized and (almost surely) strictly tighter than~\eqref{eq:cheby-ci}, when we construct infinitely many such intervals for different problems (with independent data from distributions with different means and variances), at least 95\% of those confidence intervals will cover the corresponding means. 

Our e-value based tests~\eqref{eq:rej-E-rand}, like our more powerful variant of universal inference~\eqref{eq:UI-U}, also have the same frequency interpretation as the nonrandomized tests~\eqref{eq:rej-E}. There is a simple way to interpret our use of randomization. When $\alpha=0.05$, nonrandom thresholding rules reject the null when the e-value exceeds $1/\alpha=20$. Meaning that if the e-value equals 10, we do not reject, while if it equals 20, we do. In our randomized setting, we simply view an e-value of 10 as having half the evidence as that of an e-value equaling 20, so we reject it with probability one half. Similarly, an e-value of 19.999 would not get rejected with the usual nonrandomized rules, but would get rejected with very high probability in our scheme. In other words, the rejection probability is exactly proportional to the required evidence for a definite rejection, resulting in a ``smoothed'' test, as opposed to a sharp 0-1 decision.

\subsection{The potential for p-hacking, and ideas to overcome it}

Despite the frequency interpretation discussed in the previous subsection, we recognize the potential risks for ``p-hacking'' (or hacking confidence intervals), where a naive or dishonest practitioner may reconstruct our uniformly randomized confidence interval many times (redrawing $U$), and pick the one that suits them (in order to report a narrow enough interval for their purposes),  even reporting the random seed for ``reproducibility''. This of course is not valid. Thus our intervals must be employed with care. We suggest a few points that should be kept in mind for practical applications:
\begin{itemize}
    \item If the interval constructions are coded up as part of some automated software that constructs thousands of such intervals (and perhaps acts on them) without any human involvement, then the above frequency interpretations will be preserved and~\eqref{eq:cheby-u-ci} and~\eqref{eq:hoeff-ci} yield bona fide, valid confidence intervals.
    \item In Section~\ref{subsec:truncate}, we point out that it may make sense to sometimes truncate the intervals, which increases interpretability, avoids contradicting intuition, and also reduces the extent to which p-hacking is possible.
    \item In Section~\ref{subsec:external-rand}, we point out that sometimes external randomization is not needed at all, and one can use randomness intrinsic in the data itself.
    \item Finally, we note that there may be opportunities to systematize the use of external randomness. One could create a central repository of uniform random numbers, and you can request a fixed number of them, but then have to use all of them. Independently, one could report a file with all the uniform random numbers used in an analysis, and these should pass a battery of uniformity tests (though this could itself be p-hacked, it's now a much higher bar).
\end{itemize}

\subsection{Truncation to avoid empty (or tiny) confidence intervals}\label{subsec:truncate}

   Since $\log(U)$ has its smallest possible value being $-\infty$,  the interval~\eqref{eq:hoeff-ci} could sometimes simply be the empty interval. This will only happen to (much) less than an $\alpha$ fraction of constructed intervals, since the $(1-\alpha)$ coverage property does still hold; indeed a direct calculation shows that an empty interval is constructed if and only if $\log U \leq -2\log(2/\alpha)$, which happens with probability at most $\alpha^2/4$. 
   Rare as it may be, this phenomenon may not be very useful or intuitive in practice. Thus we suggest the following alternative:
    \begin{equation}\label{eq:hoeff-ci-clt}
\bar X_n \pm  \left[\max\left( \sigma \sqrt{\frac{2\log(2/\alpha)}{n}} + \sigma \frac{\log(U)}{\sqrt{2n\log(2/\alpha)}} ~ , ~   \frac{\sigma z_{1-\alpha/2}}{\sqrt n} \right)\right],
\end{equation}
where $z_{1-a}$ is the right $a$-quantile of the standard Gaussian distribution. In short, whenever the randomized Hoeffding interval becomes smaller than the asymptotic interval based on the central limit theorem, we resort to reporting the latter. The interval in~\eqref{eq:hoeff-ci-clt} is never shorter than what the CLT reports, almost surely shorter than the original Hoeffding interval, and is nonasymptotically valid. 

Despite the fact that the Chebyshev interval is non-empty almost surely, we recognize that observing an extremely short interval, even if by chance due to randomization $U$, may also be troubling. One simple fix is to output the interval:
\begin{equation}\label{eq:cheby-u-ci-cut}
  \bar X_n \pm \frac{\sigma\max(\sqrt{U},1/2)}{\sqrt{\alpha n}}   
\end{equation}
Of course, $1/2$ can be replaced by any other constant. The above interval is still nonasymptotically valid at level $\alpha$, is almost surely tighter than Chebyshev's inequality (that is, the interval obtained from it) but never improves on it by more than a factor of 2. Indeed, the expected ratio of widths is $\EE[\max(\sqrt{U},1/2)] = 17/24 \approx 0.71$, only a mild increase from the $2/3$ value obtained earlier.

In short, enlarging the interval by truncating the random improvement may be a suitable practical middle ground.

\subsection{Using internal randomization in lieu of external randomization}\label{subsec:external-rand}

Several of the bounds in this paper were formulated in terms of external randomization $U$. However, we note that in some situations,  we can avoid the use of $U$ entirely, while maintaining the gist of the original statements. 

To elaborate, recall that many statistics in the paper, like the sample mean $\bar X_n$, are only functions of the order statistics  of the data $\mathbf X := \{X_{(1)},\dots,X_{(n)}\}$ (equivalently, of the unordered bag of data, or of the empirical distribution). For a real $x$ and 
finite set $S$ of reals, define \[\text{rank}(x; S):= \sum_{i \in S}1(x_i\leq x) / |S|.\]
Now, note that if $X_1,\dots,X_n$ are i.i.d.\ from a continuous univariate distribution, then
\[
\mathbf X \perp \text{rank}(X_n; \mathbf X).
\]
Further, the aforementioned rank  is 
uniformly distributed
on the discrete set $\{1/n,2/n,\dots,1\}$. Thus, the rank 
stochastically dominates $U$ and can therefore be used in its place in our earlier bounds. In essence, the leftover information in the data ordering, like the rank of $X_n$ within the set $\mathbf X$, can be utilized ``for free'' without affecting the distribution of the underlying main statistic (like $\bar X_n$). 

Thus, to use two examples whose expressions were recalled in the previous subsection, both
\begin{equation}\label{eq:hoeff-ci-without-u}
\bar X_n \pm  \left[ \sigma \sqrt{\frac{2\log(2/\alpha)}{n}} + \sigma \frac{\log(\text{rank}(X_n; \mathbf X))}{\sqrt{2n\log(2/\alpha)}}  \right],
\end{equation}
and
\begin{equation}\label{eq:cheby-u-ci-without-u}
  \bar X_n \pm \frac{\sigma \sqrt{\text{rank}(X_n; \mathbf X)}}{\sqrt{\alpha n}}   
\end{equation}
are valid $(1-\alpha)$ confidence intervals for the mean (of a $\sigma$-subGaussian distribution, or a distribution with variance at most $\sigma^2$, respectively).

The same observation also applies to universal inference, because the likelihood calculation in the split likelihood ratio statistic is only a function of the unordered set $\{X_1,\dots,X_m\}$ (where $m$ is the size of the first split), and hence the rank of $X_m$ within that set can be used in place of $U$ in~\eqref{eq:rand-univ-inf}.

Perhaps such uses of the uniformly-randomized Markov's inequality, which are entirely ``intrinsic'' to the data itself, may be more palatable to those who have reasons to not prefer to use ``extrinsic'' randomization by using $U$.

\section{Summary}\label{sec:conclusion}

Despite our applications being focused on
universal inference and ``betting-based'' inference, the use of the
exchangeable Markov inequality was really enabled by two properties:
\begin{itemize}
    \item The original underlying problem statement has a certain symmetry, for example the data are i.i.d.\  or exchangeable.
    \item The original method did not respect the above symmetry, by employing sample splitting, or by processing the data one at a time in a random order.
\end{itemize}
To be clear, there were some benefits to deviating from the original symmetric problem statement: in the case of universal inference, it enabled constructing a test without regularity conditions, and in the case of testing bounded means, it enabled a powerful test that was adaptive to the underlying unknown variance of the data. The loss of the problem symmetry could be regained by ``algorithmic derandomization'', that is repeating the same procedure and averaging the resulting e-values. It is in this latter step that the exchangeable Markov inequality kicks in and delivers more power to the final test.

The above bullet points apply to several other problems, for example~\cite{shekhar2021nonparametric} design betting-based tests for nonparametric two-sample testing, and the same techniques would apply to that problem as well. Similarly,~\cite{waudby2020estimating} derived the only known closed-form empirical Bernstein inequality that converges in width  \emph{exactly} to Bernstein's inequality, both of which can be improved with our uniform randomization technique, while the former can also be improved with the exchangeable Markov inequality. We omit the details for brevity.

The large improvements delivered by the uniformly-randomized Markov's inequality may be unsettling to some readers, which is why we presented a version that only uses the data itself for randomization in the previous section. We anticipate more applications and discussions about when such techniques may be appropriate (or not) to emerge with time.

\subsection*{Acknowledgments}
The authors thank Johannes Ruf, Martin Larsson, Ilmun Kim, Arun Kumar Kuchibhotla, Sivaraman Balakrishnan, Jing Lei, Neil Xu and Ian Waudby-Smith for helpful conversations and suggestions. We would 
also like to thank Robin Dunn for performing the simulation study
in Appendix~\ref{app:more_ui}.

\bibliography{exch-markov}


\appendix

\section{Further details on the relationship to \cite{huber2019halving}}\label{sec:huber}

We have seen that the additively randomized Markov's inequality in Proposition~\ref{prop:ami} is equivalent to the (multiplicatively) uniformly randomized Markov's inequality in Theorem~\ref{thm:umi}, in the sense that they can be used to derive each other. We also saw that the result~\eqref{eq:huber} by \cite{huber2019halving} makes a different claim, that is neither stronger nor weaker than Markov's inequality. 

Despite the above facts, we can show that all three results are mathematically equivalent, meaning that they can all be used to derive each other. 
To see this,
first let us recall~\eqref{eq:huber} below for simplicity: 
\[
P(X+B \geq \epsilon) \leq \EE[X]/(2\epsilon),
\]
where $B \sim U[-\epsilon,\epsilon]$. Now, rewrite the left-hand side as 
\begin{equation}\label{eq:huber-ami}
P(X \geq \epsilon-B) = P(X \geq 2\epsilon - (\epsilon+B)) = P(X \geq 2 \epsilon - 2A) = P(X/2 \geq \epsilon - A),
\end{equation}
where we used the fact that $\epsilon+B$ is distributed as $U[0,2\epsilon]$, which has the same distribution as $2A$, where $A=U[0,\epsilon]$ was defined in Proposition~\ref{prop:ami}. 

We can either apply Huber's result to the left-hand side of~\eqref{eq:huber-ami} or Proposition~\ref{prop:ami} to the right-hand side of~\eqref{eq:huber-ami} to see that the two results imply each other. Nevertheless, the ``take-home message'' behind these inequalities is quite different. Indeed, as suggested by the title of his paper, Huber's focus is on halving the bounds of Markov's inequality with his two-sided additive randomization $B$ (with implications for shape-constrained settings), while ours focus has been on improving Markov's inequality and the statistical implications of such an improvement.

\section{Proof of the time-reversed Ville inequality (Theorem~\ref{thm:reverse_ville})}
\label{app:reverse}
In this appendix, we 
provide a detailed proof of the time-reversed
Ville inequality (Theorem~\ref{thm:reverse_ville}), which is at the heart
of the proofs of the EMI (Theorem~\ref{thm:emi})
and EUMI (Theorem~\ref{thm:umi+emi}).
As previously mentioned, proofs of this result
appear under varying assumptions
in the works of~\cite{doob1940}, \cite{lee1990},
\cite{christofides1990}, and~\cite{manole2021sequential}. 
In what follows, we provide
a new self-contained proof
which easily lends itself 
to deriving the EUMI. 

Let $m \geq 1$. 
Notice first that the process $(Y_t)_{t=1}^m$ defined by
$Y_t = X_{m-t+1}$ is a forward
submartingale with respect to the 
forward filtration $\calG_t = \calE_{m-t+1}$,
$1 \leq t \leq m$. 
Indeed, $(Y_t)$ is 
adapted to $(\calG_t)$, and 
\begin{align*}
\bbE[Y_{t+1}|\calG_t]
=  \bbE[X_{m-t} | \calE_{m-t+1}] 
\geq X_{m-t+1} = Y_t,
\end{align*}
for all $t=1, \dots, m-1$.
It must then also follow that $(Y_t)_{t=1}^m$
is a forward submartingale with respect
to the filtration
$$\calF_t = \sigma(Y_1, \dots, Y_t) = 
\sigma(X_{m-t+1}, \dots, X_m),\qquad 1 \leq t \leq m.$$
With these preliminaries
in place, we turn to proving the claimed
inequality.
Given $m \geq 1$, 
define 
$\tau:= \sup\{1 \leq t \leq m: X_t \geq 1/a\}$, where $\sup\emptyset=-\infty$.  
Markov's inequality implies
\begin{equation}\label{eq:proof-ville2}
    P(\tau \geq 1) = P(X_{\tau \vee 1} \geq 1/a) \leq a \cdot \EE[X_{\tau \vee 1}].
\end{equation}
Now define  
$\eta:= m-\tau+1$, 
so that $X_{\tau \vee 1} = Y_{\eta \wedge m}$.
For all $t=1, \dots, m$, we have 
$$\{m-\tau+1 = t\} = \{\tau = m-t+1\} \in \calF_t$$
thus $\eta$ and $\eta \wedge m$ are 
stopping times with respect to $(\calF_t)$.
Furthermore, the process $(Y_t)_{t=1}^m$
is trivially uniformly integrable, hence
by Doob's optional stopping theorem for submartingales,
$$\bbE[X_{\tau \vee 1}] = \bbE[Y_{\eta\wedge m}] 
\leq \bbE[Y_m]=\bbE[X_1].$$
Returning to equation~\eqref{eq:proof-ville2}
and noting that $\{\tau \geq 1\} = \{\sup_{1 \leq t \leq m} X_t \geq 1/a\}$, we have thus shown 
$$P\left(\sup_{1\leq t \leq m} X_t \geq 1/a \right) \leq a \cdot \EE[X_1].$$ 
By the bounded convergence theorem, sending $m\to\infty$ yields our claim.
\qed

\section{A general randomized tail bound}
\label{app:bentkus}
We state and prove a 
simple randomized
tail bound which contains Theorems~\ref{thm:umi},~\ref{thm:rand-chebyshev} and~\ref{thm:rand-chernoff} as special cases, and
can be used
to derive randomized
variants of other tail bounds in the literature.
\begin{proposition} 
\label{prop:general_randomized}
Let  $X$ be a  random
variable taking values in a set $\calX \subseteq \bbR$, and let $I \subseteq \bbR_+$ be an interval. Let $f:\calX \to I$ and $g : I \to \calX$ be nondecreasing Borel-measurable
functions such that $f(g(z)) \geq z$
for any $z \in I$. Then, given a
random variable $U \sim \text{Unif}(0,1)$
independent of $X$, and  $x > 0$, it holds that
$$\bbP(X \geq g(U f(x))) \leq \frac{\bbE[f(X)]}{f(x)}.$$
\end{proposition}
The proof is exactly as before. Notice first
that $Uf(x)$ takes values in $I$ since
$U$ is supported in $[0,1]$, thus the quantity
$g(Uf(x))$ is well-defined. Furthermore, 
\begin{align*} 
P(X \geq g(Uf(x))) 
 &= P(f(X) \geq f(g(Uf(x))))  \\ 
 &\leq  P(f(X) \geq Uf(x)) \\
 &= \bbE[ P(U \leq f(X)/f(x)|X)] \leq  \bbE f(X)/f(x), 
\end{align*}
which proves the claim. \qed 


As an example, we next use 
Proposition~\ref{prop:general_randomized} to derive
uniformly-randomized analogues of Cantelli's inequality, Bernstein's inequality, and of
the empirical-Bernstein inequality.

\subsection{Uniformly-randomized Cantelli inequality}
\label{sec:cantelli}

We begin by deriving a uniformly-randomized analogue of Cantelli's inequality~\citep{cantelli1929sui},
which is a one-sided version of Chebyshev's inequality\footnote{Despite its name, Cantelli's inequality apparently originated in Chebyshev's much earlier work~\citep{tchebichef1874valeurs}; see~\cite{ghosh2002probability}. In fact, it is commonly accepted that Markov's inequality itself had already been proven by Chebyshev, who was Markov's doctoral advisor.}.
In contrast to~\eqref{eq:cheby}, it states that 
\begin{equation}\label{eq:cantelli}
    P( X - \EE X \geq k\sigma) \leq \frac1{k^2+1}.
\end{equation}
It can be improved by uniform randomization as follows:
\begin{equation}\label{eq:cantelli-u}
    P( X - \EE X \geq \sqrt{U}(k\sigma + \sigma/k) - \sigma/k) \leq \frac1{k^2+1},
\end{equation}
which we call the \emph{uniformly-randomized Cantelli inequality}.
To see that~\eqref{eq:cantelli-u} is a stronger statement than~\eqref{eq:cantelli}, rewrite the left-hand side as
\[
P( X - \EE X \geq k\sigma - (1-\sqrt{U})(k\sigma + \sigma/k)),
\]
and note that $(1-\sqrt{U})$ is positive. 

The proof of~\eqref{eq:cantelli-u} is a simple
consequence of   Proposition~\ref{prop:general_randomized}.
Let $x = \sigma k$ and $u = \sigma/k$.
Taking $f(y) = (y+u)^2$ for all $y \in \bbR$ and 
and $g(z) = \sqrt z - 1$ for all $z \in \bbR_+$, we have
\begin{align*}
P\left(X-\bbE X \geq \sqrt U (k\sigma + \sigma/k) - \sigma/k\right)
 &= P\big(X-\bbE X\geq g(Uf(x))\big) \\
 &\leq \frac{\bbE[f(X-\bbE X)]}{f(x)} \\
 &= \frac{\sigma^2+u^2}{(x+u)^2} 
  = \frac{1 + 1/k^2}{2 + k^2 + 1/k^2} = \frac 1 {1+k^2},
\end{align*} 
as claimed. 


\subsection{Uniformly-randomized Bernstein inequality}
\label{sec:bernstein}

We next derive a uniformly-randomized
analogue of the classical Bernstein inequality 
(see for instance Proposition 2.10 of~\cite{wainwright2019}). 
We will assume that the random variable~$X$
satisfies the Bernstein condition, namely
\begin{equation}
\label{eq:bernstein_condition}
\big| \bbE(X-\bbE X)^k \big| \leq \frac 1 2 k! \sigma^2 b^{k-2}, \quad \text{for } k=2,3,\dots
\end{equation}
for some $\sigma,b > 0$.
In particular, if $X$ is a bounded
random variable with variance at most 
$\sigma^2$ and satisfying $|X - \bbE X| \leq 1$, then the above condition
holds for $b=1/3$.
Bernstein's inequality states that, 
under condition~\eqref{eq:bernstein_condition},
for all $\epsilon > 0$,
\begin{equation}
\label{eq:bernstein}
P\left( X - \bbE X \geq \epsilon\right)
\leq \exp\left(-\frac{\epsilon^2}{2(\sigma^2 + \epsilon b)}\right).
\end{equation}
  Setting the right-hand side
  equal to $\alpha \in (0,1)$, 
the above implies
$$P\left(X - \bbE X \geq
\sqrt{2\sigma^2 \log(1/\alpha)} + 
2 b \log(1/\alpha)\right) \leq \alpha.$$
Consequently, given $n$ i.i.d.\ samples of $X$, the following is a $(1-\alpha)$-confidence interval for $\EE X$:
\begin{equation}\label{eq:bern-ci}
    \bar X_n \pm \left[\sqrt{\frac{2\sigma^2 \log(2/\alpha)}{n}} + 
\frac{2 b \log(2/\alpha)}{n}\right].
\end{equation}
It is well known that the i.i.d.\ assumption can be relaxed into a martingale dependence assumption requiring neither the independence aspect nor the identically distributed aspect, but we omit this generalization for simplicity.

Our {\it uniformly-randomized Bernstein inequality} reads
\begin{equation}
\label{eq:randomized-bernstein}
P\left(X - \bbE X \geq \left(\sigma^2 + b\sqrt{\sigma^2\log(1/\alpha)} 2b^2\log(1/\alpha)\right) \log U
  + \sqrt{2\sigma^2\log(1/\alpha)} + 2b \log(1/\alpha)\right) \leq \alpha, 
\end{equation}
where $U$ is a random variable that is independent of $X$ and (stochastically larger than) uniform on $[0,1]$. 
Since $\log U < 0$ almost surely, the above inequality provides a strict and almost sure improvement of Bernstein's inequality (and it recovers Bernstein's inequality by substituting $U=1$). Recall also that $\bbE[\log U]=-1$ and $\Var[\log U]=1$, giving an idea of the extent of the improvement.

To prove~\eqref{eq:randomized-bernstein},
let $|\lambda| \leq b^{-1}$,
$f(x) = \exp(\lambda x)$ 
for all $x \in \bbR$,
and $g(z) = \log z / \lambda$ for all $z > 0$.
By Proposition~\ref{prop:general_randomized},
we have for all $x > 0$,
$$P\left(X - \mu \geq \frac{\log U}{\lambda}
  + x\right) \leq e^{-\lambda x} \bbE[e^{\lambda(X-\bbE X)}] \leq \exp \left\{ \frac{\lambda^2\sigma^2}{2(1-b|\lambda|)}-\lambda x\right\}, $$
where we used the fact that 
  $\bbE[e^{\lambda(X-\bbE X)}] \leq \exp(\lambda^2\sigma^2 / 2(1-b|\lambda|))$ for all $|\lambda| \leq b^{-1}$
  under the Bernstein condition (cf. Proposition 2.10
  of~\cite{wainwright2019}). 
  Now, letting $\lambda = (bx+\sigma^2)^{-1}$ 
  and simplifying the above expression, we obtain
 \begin{equation}
 P\left(X - \mu \geq (bx+\sigma^2) \log U
  + x\right) \leq 
\exp\left\{-\frac{x^2}{2(\sigma^2+bx)}\right\},
\end{equation}
which can be viewed as another form of our randomized Bernstein inequality.
Setting $x = \sqrt{2\sigma^2 \log(1/\alpha)} + 
2b\log(1/\alpha)$ leads to the claimed inequality. 

Combining the above with Lieb's inequality, one directly obtains a randomized matrix-Bernstein inequality as well; we omit the details for brevity.

\subsection{Uniformly-randomized empirical Bernstein inequality}
\label{sec:empirical_bernstein}

Bernstein's inequality is not always practically applicable due to the need to know $\sigma$. When the data are bounded, one can construct so-called empirical 
Bernstein (EB) inequalities that only depend on the data. There are several such EB inequalities in the literature, but we present below a randomized variant of a recent EB inequality by~\cite{waudby2020estimating}, because it is the only one that we are aware of whose corresponding confidence interval width exactly matches the first order term in~\eqref{eq:bern-ci}.

Going forward, suppose that the i.i.d.\ data lie in $[0,1]$; this is done for simplicity and without loss of generality. Define
\[
\psi(\lambda) := (- \log (1-\lambda) - \lambda) \quad \text { for } \lambda \in [0, 1),
\]
and the instantaneous empirical variance as
\[
v_t := (X_t - \widehat \mu_{t-1})^2, \text{ where } \widehat \mu_t := \frac{\tfrac1{2} + \sum_{i=1}^t X_i}{t + 1}
\]
The following is then a $(1-\alpha)$-CI for $\EE X$:
\begin{equation}\label{eq:eb-ci}
     \frac{\sum_{t=1}^n \lambda_t X_t}{\sum_{t=1}^n \lambda_t} \pm \left[ \frac{\log(2/\alpha) +  \sum_{t=1}^n v_t \psi(\lambda_t) }{\sum_{t=1}^n \lambda_t} \right] ,
\end{equation}
where $\lambda_t \in (0,1)$ is a function of $X_1,\dots,X_{t-1}$ that is set as follows:
\begin{equation} 
\nonumber
\lambda_t := \sqrt{\frac{2 \log (2/\alpha)}{\widehat \sigma_{t-1}^2 n }} \land \frac12, ~~~ \widehat \sigma_{t}^2 := \frac{\tfrac1{4} + \sum_{i=1}^t (X_i - \widehat \mu_i)^2}{t + 1}.
\end{equation}
The proof follows by observing that
\[
M_t := \prod_{i=1}^t \exp\left \{ \lambda_i (X_i-\mu) - v_i \psi(\lambda_i) \right \} 
\]
is a nonnegative supermartingale with initial value $M_0=1$, and thus $M_n$ is an e-value. Applying Markov's inequality to $M_n$, and rearranging, yields one side of~\eqref{eq:eb-ci}, and a union bound with $-\lambda_i$ in place of $\lambda_i$ yields the other side. Instead, applying UMI in place of Markov's inequality, we obtain the uniformly-randomized empirical Bernstein confidence interval:
\begin{equation}\label{eq:rand-eb}
     \frac{\sum_{t=1}^n \lambda_t X_t}{\sum_{t=1}^n \lambda_t} \pm \left[ \frac{\log(2/\alpha) + \log U +  \sum_{t=1}^n v_t \psi(\lambda_t) }{\sum_{t=1}^n \lambda_t} \right],
\end{equation}
which is almost surely tighter than~\eqref{eq:eb-ci} since $\log U < 0$ with probability one.

\begin{remark}
    Denoting the expression in~\eqref{eq:eb-ci} by $C_n$, it turns out that $\bigcap_{i \leq n} C_i$ is also a valid $(1-\alpha)$-confidence interval. This is obtained by applying Ville's inequality to the supermartingale $M$, in place of Markov's inequality. It may be a priori unclear whether $\bigcap_{i \leq n} C_i$ is tighter than~\eqref{eq:rand-eb} or not, but Figure~\ref{fig:betting_umi} suggests a clear win for UMI over Ville.
\end{remark}

\section{Additional simulation results}
\label{app:more_sims}
In this appendix, we 
include several simulation
results which were omitted
from Section~\ref{sec:simulations}. 

\subsection{Additional simulation results
from Section~\ref{sec:sim_evals}}

\begin{figure}[H]
\centering
\includegraphics[width=\textwidth]{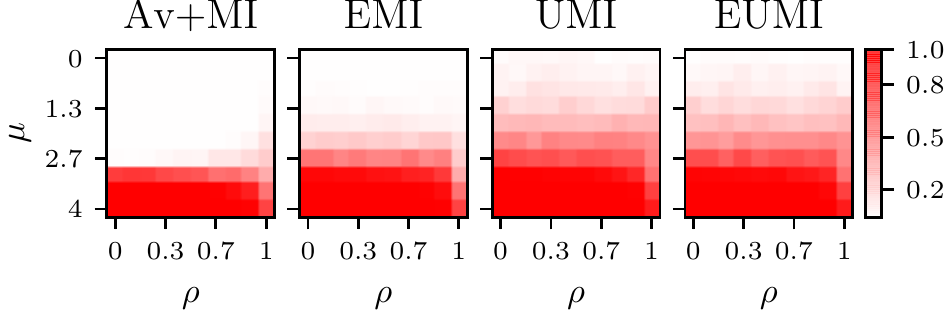}
\caption{
\label{fig:e_val_all}
Individual rejection proportions of the 
methods Av+MI, EMI, UMI, and EUMI, 
in the simulation study of Section~\ref{sec:sim_evals}. }
\end{figure}

\begin{figure}[H]
\centering
\includegraphics[width=\textwidth]{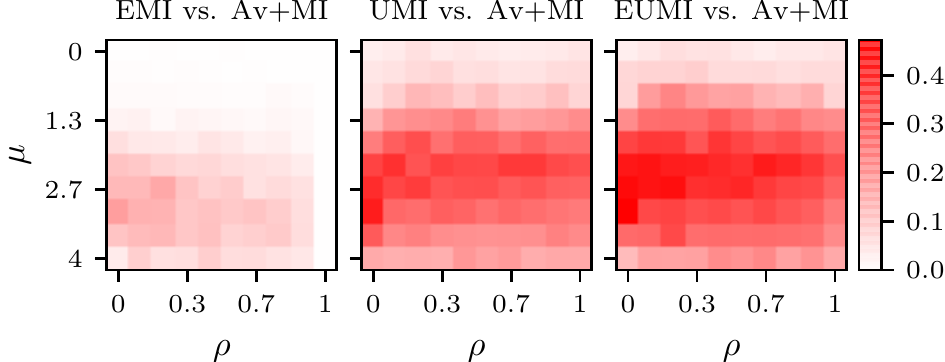}
\caption{\label{fig:K2-e-val} Comparison of the
rejection proportions $\pi_{\mathrm{AvMI}}$, $\pi_{\mathrm{EMI}}$, $\pi_{\mathrm{UMI}}$,
and $\pi_{\mathrm{EUMI}}$. 
of the procedures~\eqref{eq:naive-e-avg}, \eqref{eq:u-e-avg}, ~\eqref{eq:eval-perm},
and~\eqref{eq:u-e-perm-avg}, 
for rejecting the null hypothesis~$H_0:\mu\leq 0$
based on the e-values~\eqref{eq:sim-e-vals}, but now with $K=2$
rather than $K=100$.
For varying
values of~$\rho$ and~$\mu$, the left-hand
side plot represents 
the difference $\pi_{\mathrm{EMI}} - \pi_{\mathrm{AvMI}}$, 
the middle plot represents
$\pi_{\mathrm{UMI}} - \pi_{\mathrm{AvMI}}$, 
and the right-hand side plot represents 
$\pi_{\mathrm{EUMI}} - \pi_{\mathrm{AvMI}}$.}
\end{figure}

\begin{figure}[H]
\centering
\includegraphics[width=\textwidth]{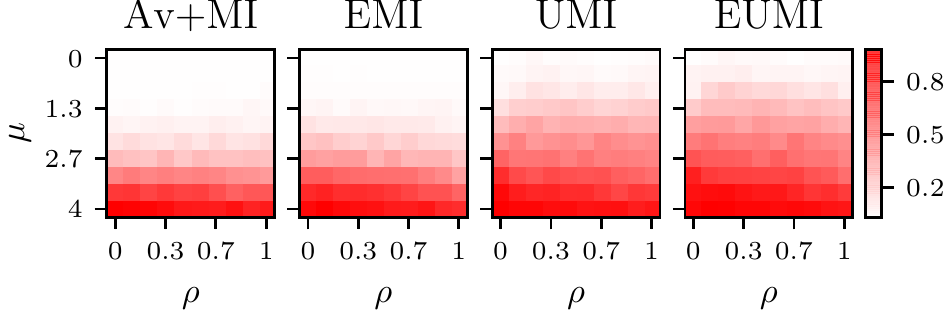}
\caption{
\label{fig:K2-e_val_all}
Individual rejection proportions of the 
methods Av+MI, EMI, UMI, and EUMI, 
in the simulation study of Section~\ref{sec:sim_evals},
but now with $K=2$ rather than
$K=100$. }
\end{figure}

\clearpage 

\subsection{Additional simulation results
from Section~\ref{sec:sim_betting}}

\begin{figure}[H]
\centering
\includegraphics[width=0.8\textwidth]{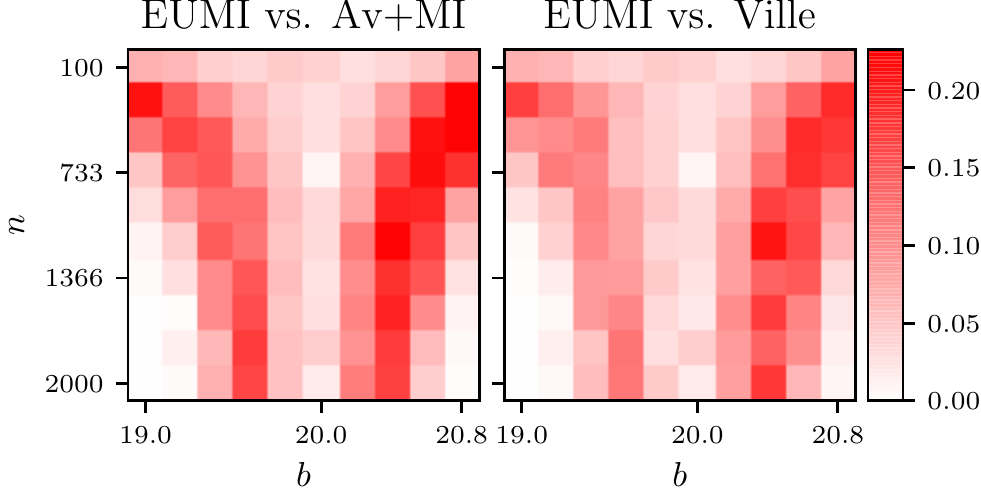}
\caption{\label{fig:betting_eumi} Comparison of the
rejection proportions $\pi_{\mathrm{AvMI}}$, $\pi_{\mathrm{Ville}}$, and $\pi_{\mathrm{EUMI}}$ 
of the respective procedures~\eqref{eq:betting_MI}, \eqref{eq:bet-ville2} and~\eqref{eq:bet-umi-emi},
for rejecting the null hypothesis~$H_0:b=20$
based on the statistic $M_n$.
For varying
values of~$b$ and~$n$, the left-hand
side plot represents 
the difference $\pi_{\mathrm{EUMI}} - \pi_{\mathrm{AvMI}}$, 
and the right-hand side plot represents 
the difference $\pi_{\mathrm{EUMI}} - \pi_{\mathrm{Ville}}$.}
\end{figure}

\begin{figure}[H]
\centering
\includegraphics[width=0.8\textwidth]{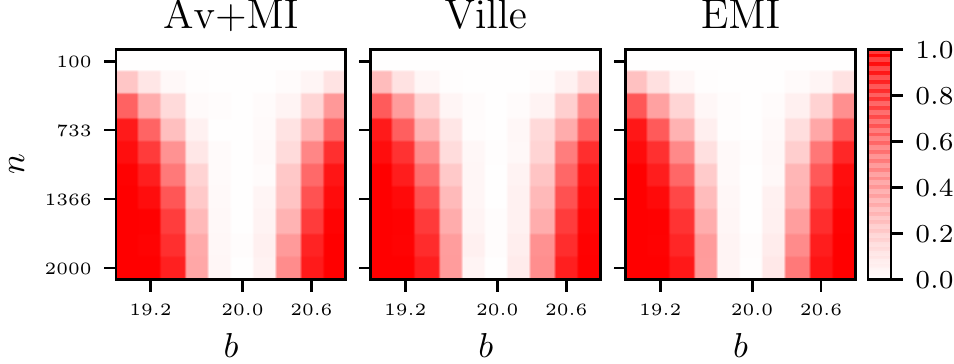}
\includegraphics[width=0.5\textwidth]{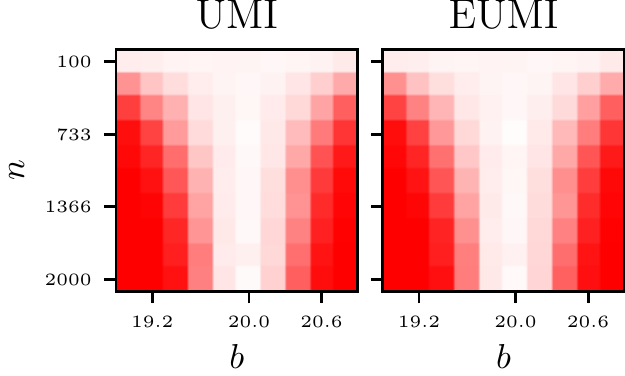}
\caption{\label{fig:betting_all}
Individual rejection
proportions of the methods Av+MI, Ville, EMI, UMI, and EUMI, in the simulation study of Section~\ref{sec:sim_betting}.}
\end{figure}

\begin{figure}[H]
\centering
\includegraphics[width=0.5
\textwidth]{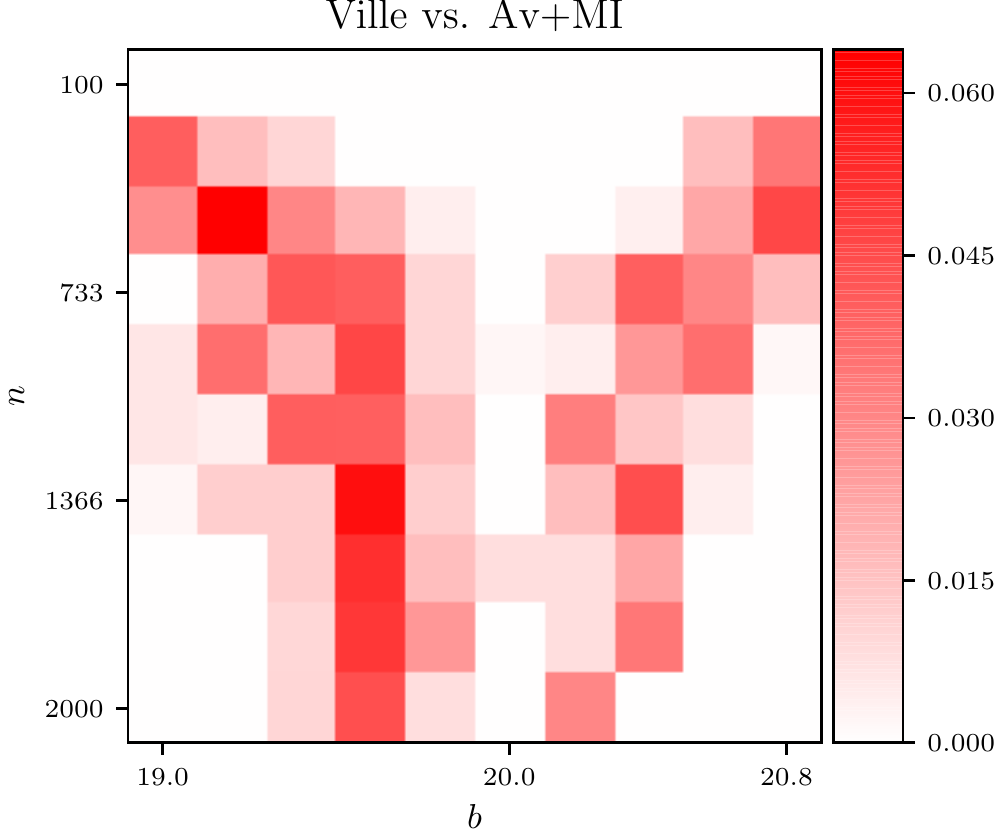}
\caption{\label{fig:betting_MI_Ville}
Plot of $\pi_{\mathrm{Ville}} - \pi_{\mathrm{AvMI}}$, where
$\pi_{\mathrm{Ville}}$ and 
$\pi_{\mathrm{AvMI}}$ respectively
denote the rejection proportions
of the methods Ville and MI
in the context of the simulation study
in Section~\ref{sec:sim_betting}. }

\end{figure}

\subsection{Simulation
results for Gaussian Universal Inference}
\label{app:more_ui}
In this appendix, we briefly
report further simulation results
for the Universal
Inference method. Let
$X_1, \dots, X_n$
be an i.i.d. sample 
from the  $\calN(\theta^*,I_2)$
distribution, where the unknown
parameter $\theta^*\in \bbR^2$ is fixed to $\theta^* = 0$. 
Figure 1 of~\cite{dunn2021gaussian}
reports six simulated
confidence sets for $\theta^*$
based on the split likelihood
ratio statistic~\eqref{eq:split-lrt}, the closely-related cross-fit likelihood
ratio statistic~\citep{wasserman2020universal}, the subsampling likelihood ratio
statistic (appearing
in~\eqref{eq:UI}), and the classical likelihood
ratio (LRT) statistic.  
We reproduce this figure for the reader's convenience as Figure~\ref{fig:gaussian_ui_orig} below. Six simulations
are reported for a single
random sample; the differences
across simulations arise merely
from the randomness inherent in the
various methods.
\begin{figure}[H]
\begin{center} 
\includegraphics[width=0.9\textwidth]{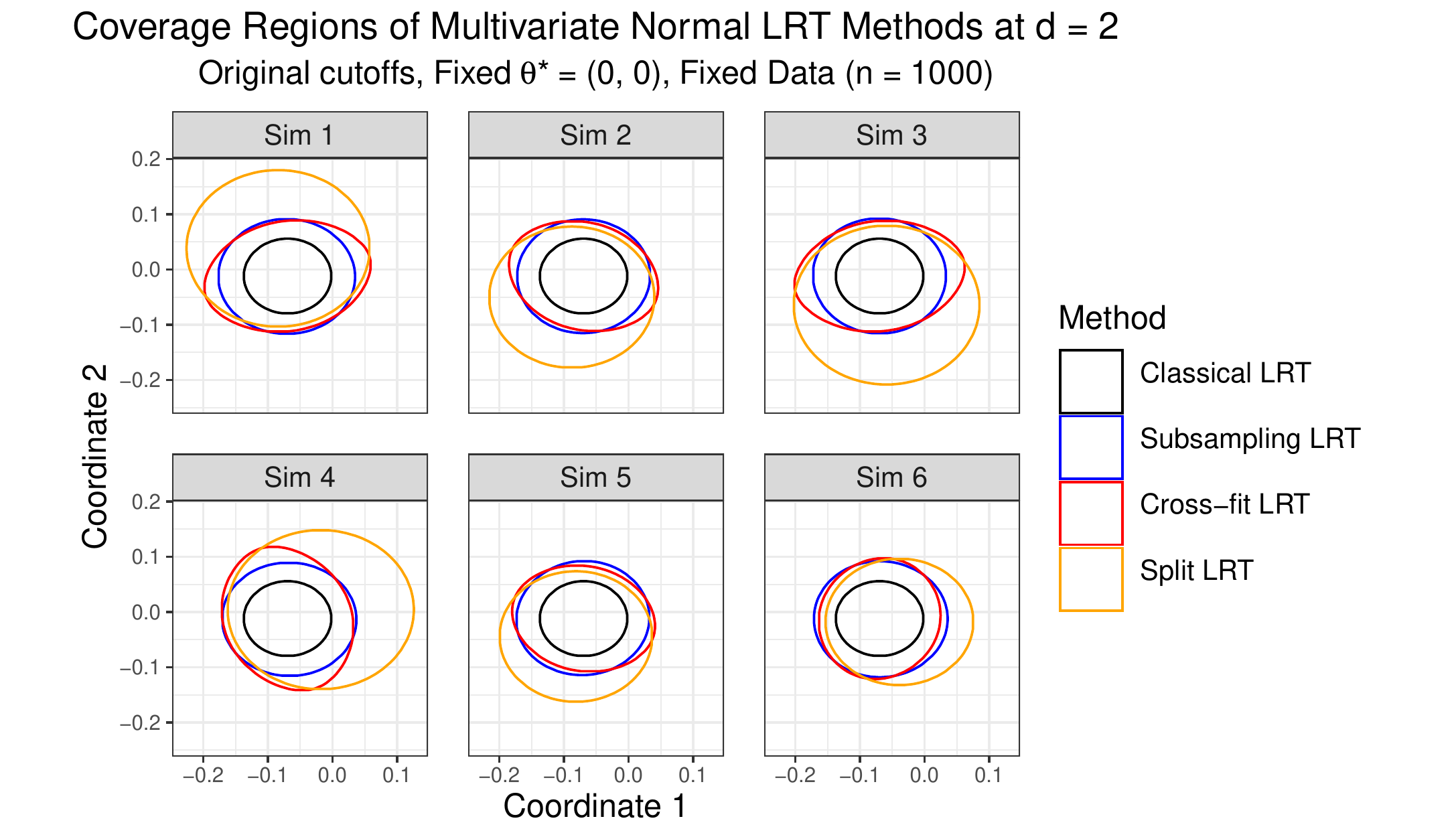}
    \caption{Confidence sets for $\theta^*$ based on a single sample of size $n=1,000$, 
    as described in~\cite{dunn2021gaussian}.}
    \label{fig:gaussian_ui_orig}
\end{center}
\end{figure}
In Figure~\ref{fig:gaussian_ui_umi}, we report
the confidence sets produced by the same four methods, but 
now replacing the MI by  the UMI 
in the definitions
of three of them (excluding the classical LRT). 
Each of these
three confidence sets
is visibly smaller when using
the UMI rather than the MI.
In some cases, the difference
is striking, with the volume 
of the UMI-based confidence sets
being nearly equal
to that of the
LRT confidence set.

Of course, these are just 6 runs (and we did not run it more than 6 times), but one still gets a sense of the achieved improvement. If one gets lucky (meaning draws a small value of $U$), the confidence sets might even get smaller than the LRT's. As discussed in Section~\ref{sec:sum}, if one would like to ensure (for the sake of intuition) that the methods based on UMI never \emph{beat} the LRT, one could always take the union of the achieved sets with that of the LRT.

\begin{figure}[H]
\begin{center}
    \centering
\includegraphics[width=0.9\textwidth]{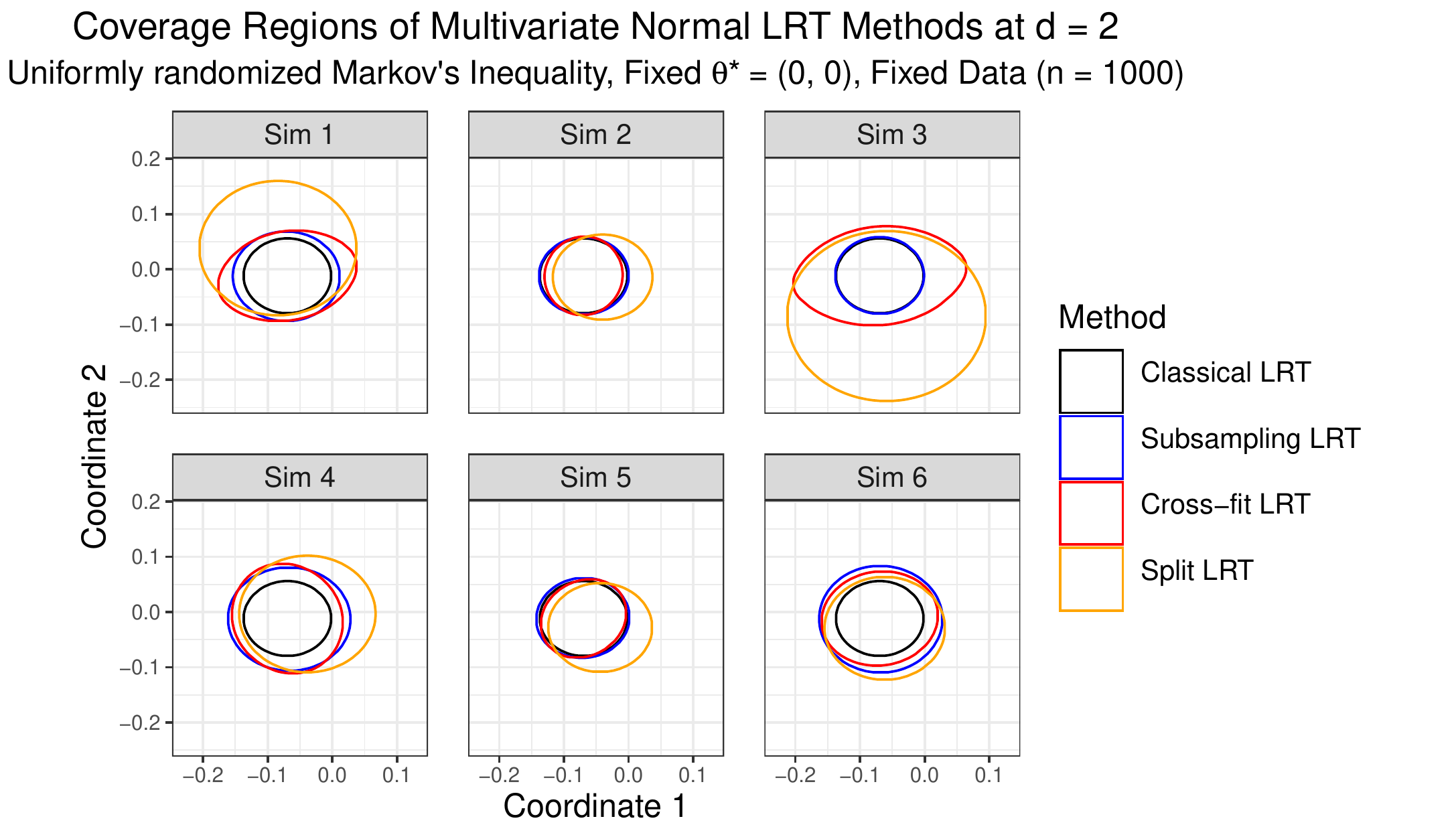}
    \caption{Confidence sets for $\theta^*$ based on 
    the same
    sample of size $n=1,000$
    as in Figure~\ref{fig:gaussian_ui_orig}, 
    but now using the UMI
    rather than the MI. The Universal Inference
    confidence sets are
    uniformly tighter than those presented
    in Figure~\ref{fig:gaussian_ui_orig}, 
    in some cases strikingly so.}
    \label{fig:gaussian_ui_umi}
\end{center}
\end{figure}



\end{document}